\documentclass[11pt,amstex]{article}
\usepackage{amsmath}
\usepackage{amssymb}
\usepackage{amsfonts}
\usepackage{amscd}
\usepackage[usenames]{color}
\usepackage{setspace}
\usepackage{verbatim}
\usepackage{stmaryrd}

\input xypic
\xyoption{all}

\def\a{\alpha}
\def\b{\beta}

\def\de{\delta}
\def\d{\delta}
\def\ep{\varepsilon}
\def\la{\lambda}
\def\si{\sigma}

\def\Ka{K_\alpha}

\def\be{\begin{equation}}
\def\ee{\end{equation}}
\def\bear{\begin{eqnarray}}
\def\eear{\end{eqnarray}}
\def\best{\begin{eqnarray*}}
\def\eest{\end{eqnarray*}}
\def\pf{{\bf Proof}: }
\renewcommand{\theequation}{\arabic{section}.\arabic{equation}}
\newtheorem{theorem}{Theorem}[section]
\newtheorem{prop}[theorem]{Proposition}

\newtheorem{lemma}[theorem]{Lemma}

\newtheorem{defn}[theorem]{Definition}

\newtheorem{Vthm}[theorem]{Vanishing Theorem}

\newtheorem{Theorem A}[theorem]{Theorem A}
\newtheorem{remark}[theorem]{Remark}
\newenvironment{rem}{\begin{remark}\rm}{\end{remark}}
\newtheorem{example}[theorem]{Example}
\newenvironment{ex}{\begin{example}\rm}{\end{example}}

\def\non{\noindent}
\def\pf{\non {\bf Proof. }}
\def\qed{\nopagebreak \hskip .1in { $\Box$ }\penalty10000 %
\hskip\parfillskip \par  }

\def\del{\overline \partial}
\def\bd{\partial}

\def\dim{\mbox{\rm dim }}

\def\ind{\mbox{Index\,}}
\def\ker{\mbox{ker\,}}
\def\cok{\mbox{coker\,}}

\def\ov#1{\overline{#1}}
\def\Z{{\mathbb Z}}
\def\R{{\mathbb R}}
\def\P{{\mathbb P}}
\def\Q{{\mathbb Q}}
\def\cx{{\mathbb C}}

\def\ev{\mbox{\rm ev}}
\def\im{\mbox{\rm Im}}

\def\Hom{{\rm Hom}}

\def\Ka{K\"{a}hler }

\def\E{{\cal E}}
\def\F{{\cal F}}

\def\H{{\cal H}}
\def\N{{\cal N}}
\def\O{{\cal O}}

\def\bU{{\overline{\cal U}}}

\def\M{{\cal M}}
\def\bM{\ov\M}

\def\ind{{\rm index\;}}

\newcommand{\CM}{\overline{{\cal M}}}
\newcommand{\CU}{\overline{{\cal U}}}

\newcommand{\cred}{\color{red}}

\def\Ob{{{\cal O}b}}

\def\Map{{\cal{M}}{\it{ap}}}

\def\vir{\mbox{\rm \scriptsize vir}}

\def\ss{\scriptscriptstyle}

\def\H{{\cal{H}}}

\def\dist{\text{\rm{dist}}}
\def\diam{\text{\rm{diam}}}
\def\piN{\pi_{\scriptscriptstyle N}}

\def\ev{\mbox{\it{ev}}}
\def\hev{\hat{ev}}
\def\CP{{\cal P}}
\def\vv{\interleave}
\def\ssetminus{\hspace{-.03cm}\setminus\hspace{-.04cm}}

\title{{\bf An obstruction bundle relating Gromov-Witten invariants of curves and K\"{a}hler surfaces}}
\author{Junho Lee and Thomas H. Parker\thanks{partially supported by
the N.S.F.}}

\date{\empty}
\addtocounter{section}{0}
\begin{document}

\maketitle

\begin{abstract}
\medskip

In \cite{LP} the authors defined symplectic  ``Local Gromov-Witten
invariants'' associated to spin curves and showed that the GW
invariants of a  K\"ahler surface $X$ with $p_g>0$ are a sum of such
local GW  invariants.   This paper  describes how the local GW
invariants arise from an obstruction bundle (in the sense of Taubes
\cite{T}) over the space of stable maps into curves. Together with the results of \cite{LP},  this reduces the calculation  
of the GW invariants of complex surfaces to computations in the  GW theory of curves.
\end{abstract}
\vspace{.8cm}

\setcounter{section}{0}

For a compact  K\"{a}hler surface $X$, a holomorphic 2-form $\a$ is
a section of  the canonical bundle whose zero locus is a canonical
divisor $D$.    Several years ago, the first author observed that
each such  2-form $\a$ naturally induces an almost complex structure
$J_\a$  that satisfies a remarkable property:

\medskip

\noindent {\bf Image Localization Property:}  {\em If a $J_\alpha$-holomorphic map represents a (1,1) class, then its image lies in
$D$. }

\medskip

\noindent After further perturbing to a generic $J$ near $J_\a$, the
images of all $J$-holomorphic maps cluster  in $\ep$-neighborhoods
of the  components of the canonical divisor  $D$.   This implies
that
 the Gromov-Witten invariant of $X$  is a sum
$$
GW_{g,n}(X,A)\ =\ \sum_k GW^{loc}_{g,n}(D_k, d_k[D_k])
$$
over the connected components $D_k$ of $D$ of ``local invariants''
that count the contribution of maps whose image lies near $D_k$.

When $D$ is smooth, the local invariants depend only on   the normal
bundle $N$ to $D\subset X$.  By the adjunction formula, $N$ is a
holomorphic square root of the canonical bundle $K_D$, that is, $N$
is a  theta-characteristic of the curve $D$ and the pair
$(D,N)$ is a  spin curve. The total space $N_D$ of $N$ has a
tautological holomorphic 2-form $\a$ whose zero locus is the zero
section $D\subset N_D$. For perturbations of the corresponding
$J_\a$,  all $J$-holomorphic maps cluster around the zero section. 
These clusters define local GW invariants of the spin curve $(D,N)$
which, the authors proved in
 \cite{LP},  depend only on  the parity of $N$
(i.e. on $h^0(N)$ mod 2).  Altogether, we have
\begin{equation}\label{key}
GW_{g,n}(X,A) \ =\ \sum_k (i_k)_*\,GW^{loc}_{g,n}(N_{D_k},d_k) \vspace{.3cm}
\end{equation}
where $(i_k)_*$ is the induced map from the inclusion $D_k\subset
X$. Thus the computation of the GW invariants of  \Ka surfaces with
$p_g>0$ is reduced to the problem of calculating the local
invariants of spin curves.

Recently, Kiem and Li \cite{KL} defined
the local invariants by algebraic geometry methods and proved the formulas for degree 1 and 2 local invariants
conjectured by  Maulik and Pandharipande \cite{MP}. The first author
\cite{L2} also reproved those formulas by adapting  the symplectic sum
formula of \cite{IP2} to local GW invariants.

Because (\ref{key}) applies to all GW invariants, not
just those of the ``embedded genus'', one cannot apply
Seiberg-Witten theory. Nor can they be computed by the usual methods
of algebraic geometry, such as localization and
Grothendieck-Riemann-Roch,  because {\em the linearized
$J_\a$-holomorphic map equation is not complex linear}.  In
particular, when $\mbox{genus}(D_k)>0$  the local invariants in (\ref{key}) are not the same as
the ``local GW invariants'' used to study Calabi-Yau 3-folds
\cite{BP} or the ``twisted GW invariants'' defined by Givental.

While the local invariants are defined in terms of the GW theory of
the (complex) {\em surface} $N_D$ one would like, as a step toward
computation,
 to recast them in terms of  the much better-understood GW theory of {\em curves} (cf. \cite{OP}).   This paper uses  geometric  analysis arguments to prove that the local GW invariants  of a spin curve $(D,N)$ arise from a cycle in the space of stable maps into the curve $D$.  The cycle is defined by constructing  an ``obstruction bundle''.   While the basic idea is clear and intuitive, the construction is difficult because of technical issues involving the construction of a complete space of maps.

The intuition goes like this.  The tautological 2-form $\a$ on $N_D$
determines an almost complex structure $J_a$.   By the  Image
Localization Property,  the space   of stable $J_\alpha$-holomorphic
maps  into $N_D$ representing $ d[D]$ is the same as   the space
of  degree $d$ stable maps into $D$:
\begin{equation}
\label{M2} \bM^{J_\alpha}_{g,n}(N_D, d[D])\ =\ \bM_{g,n}(D, d).
\end{equation}
Counting dimensions, one sees that  the formal dimension of
$\bM_{g,n}(D, d)$  is exactly twice the dimension of the virtual
fundamental class that defines the local GW invariants of $N_D$ (when $n=0$). The
dimensions do not match because $J_\a$ is not generic.  Perturbing
$J_\a$ to  a generic $J$ effectively reduces the space of maps to  a
half-dimensional cycle in $\bM_{g,n}(D, d)$ that defines the local
GW invariants of the spin curve $(D,N)$.  To understand this
reduction, we use
  another remarkable property of the   $J_\alpha$-holomorphic maps:

\medskip

\noindent {\bf Injectivity Property:}  {\em The linearization of the
map  $f\mapsto \del_{J_\alpha}f$, when restricted to the normal
bundle,  is an elliptic operator  $L_f$  whose kernel vanishes for
every $J_\alpha$-holomorphic map ($\alpha\neq 0$).}

\medskip

The injectivity property implies that the vector spaces $\cok L_f$
have constant dimension. The construction of  Section~6     shows
that these cokernels  form a locally trivial
vector  bundle
\begin{equation}
\label{Ob}
\begin{CD}
\Ob\\
@VVV\\
\bM_{g,n}(D, d)
\end{CD}
\end{equation}
over the space (\ref{M2}) whose rank is the formal GW dimension for
surface. This is an  ``obstruction' bundle'' in the sense of Taubes:
the Implicit Function Theorem shows that  the space of perturbed
holomorphic maps is diffeomorphic to the subset of $\bM_{g,n}(D, d)$
given by the zero set of a certain section of $\Ob$.  As always, we have the associated map
\begin{equation}
\label{1.evmap}
\hev=st\times \ev:\bM_{g,n}(D, d) \to \bM_{g,n}\times D^n
\end{equation}
whose first factor is the stabilization map  and whose second factor
records the images of the marked points.    In this context, the GW invariants of the curve $D$ are defined
 as the image of the virtual  fundamental class $[\CM_{g,n}(D,d)]^{ \vir}$ under $\hev_*$.  Our main result is that the
 local GW invariants are defined from these by capping with the Euler class of this obstruction bundle:

\bigskip\medskip

\noindent{\bf Main Theorem}
 \label{theoremB}
{\em There  is an
oriented real bundle $\Ob$ over $\CM_{g,n}(D,d)$  isomorphism class
of $\Ob$ depends only on the parity of the spin curve $(D,N)$.  When $\mbox{genus}(D)>0$
\begin{equation}\label{lp3}
 GW_{g,n}^{loc}(N_D,d) \,=\,  \hev_*\big(\,[\CM_{g,n}(D,d)]^{ \vir} \cap e(\Ob)\,\big)
\end{equation}
where  the virtual fundamental class $[\CM_{g,n}(D,d)]^{ \vir}$
defines the GW invariants of the curve $D$.   }

\medskip

This formula describes how the local invariants evaluate on elements of  the cohomology of $\bM_{g,n}\times D^n$.  But in the context of this theorem, all descendant classes are pullbacks of  classes in $H^*(\bM_{g,n}\times D^n)$ (see Section~1).  Thus the equality (\ref{lp3}) applies to descendant classes.

This theorem is a step  toward computing the GW invariants of  
 minimal \Ka surfaces with $p_g>0$.   For non-minimal surfaces, one  would also need  a version of the Main Theorem for the  local GW  
invariants of an exceptional curve:  when $N$ is the bundle $\O(-1)$ over $D=\P^1$.   This case is fundamentally different: it is simpler because the Image Localization and Injectivity Properties hold for  the standard complex structure,  but is more complicated  because a  lemma essential  for the analysis -- Lemma~\ref{lemma-1} below -- does not apply when $D=\P^1$.  As a result, the right-hand side of (\ref{lp3}) is well-defined and the associated GW invariants are computable \cite{FP}, but our analysis does not  show the equality in (\ref{lp3}).   This $D=\P^1$ case will be analyzed elsewhere.

\bigskip

The main theorem would not be difficult  to prove if  $\Ob\to \CM_{g,n}(D,d)$ were a smooth vector bundle over a manifold, or if one could perturb to the smooth situation while retaining the image localization and injectivity properties.   Unfortunately, there do not currently exist theorems or techniques for smoothing the compactified moduli space.   As a result, the proof of the main theorem must  address two significant technical issues: (i)  the lack of information about the structure  of the moduli space near its boundary, and (ii) the local triviality of the obstruction bundle.

 To deal with issue (i), we first  perturb the  $J_\alpha$-holomorphic map equation in the manner of Ruan-Tian \cite{RT2}.  The resulting 
 moduli space $\CM=\CM_{g,n}(D,d)$ of $(J,\nu)$-holomorphic maps into $D$ consists of a  smooth top stratum $\M$ and   boundary strata
of codimension at least two.  We fix a small neighborhood $U$ of the boundary and consider its complement 
 $\M_U=\M\ssetminus U$.   Then $\M_U$ is a compact smooth oriented manifold with boundary, so defines a {\em relative} homology class in $\Map$.   In Section~3 we prove that this ``relative  virtual fundamental class'' defines the same Gromov-Witten invariants as other standard definitions, including the one used in algebraic geometry.
 
 On the compact manifold $\M_U$,  the $(J,\nu)$-holomorphic map equation defines a section of the obstruction bundle, and 
  the Euler class $e(\Ob)$ is Poincar\'{e} dual to the zero set of any  section that is transverse to zero.  We achieve transversality by adding 
 a second Ruan-Tian perturbation term.  In  Section~7 we prove a generalized image localization theorem  that shows that this second perturbation leaves $\M$ unchanged.  The main theorem is proved in  Section~8 by showing that the zero set  of this section defines a rational homology class
 that  is equal to the local GW invariants by cobordism and on the other hand is equal to the right-hand side of (\ref{lp3}) by  Poincar\'{e} duality for $\M_U$.

The analysis aspects of the proof are aimed at issue (ii) above:  proving that the obstruction bundle is locally trivial.   The key difficulty is that the normal component of the linearization of the $J_\alpha$-holomorphic map equation at a map $f$  is an operator of the form $L_f=\del +A(df)$, and $df$ is not pointwise bounded in  
the topology of the Gromov compactness theorem.  Thus it is not clear  whether  $L_f$ is continuous in $f$ on the space of smooth  maps -- a fact we need in order to show that $\Ob$, which is essentially $\cok L_f$, is locally trivial  (in the literature, continuity is often implicitly assumed). For this purpose, we introduce a stronger topology on the space of maps in Section~2  and prove a strengthened version of Gromov compactness theorem.   Sections~4 and 5 develop the needed analysis results to show that $L_f$ is continuous in $f$ as an operator on appropriate weighted Sobolev spaces.   These results are used in Section~6 to define the obstruction bundle and prove that it is locally trivial.

\medskip

For computations, one would like to express the obstruction bundle $\Ob$   in terms of algebraic
geometry.  At the K\"{a}hler  structure $J_0$ on $N_D$, 
 $L_f$  is the $\del$-operator on the bundle
$f^*N$, the fiber  $\Ob_f$ is  $H^{0,1}(f^*N)$, and the injectivity property suggests that $h^0(f^*N)=0$.  This would imply that $\Ob$ is  the index bundle $\mbox{ind\,}\del$.  However, when  $D$ has genus $h>1$ the  injectivity property does not hold for $J_0$  and, as in Brill-Noether theory,   $h^{0,1}(f^*N)$ can jump up at special maps. (Pandharipande and Maulik showed us a specific example where such a jump necessarily occurs in the moduli space, and a similar example appears in \cite{KL}.)  Thus  $\Ob$ is not in general equal to  the $\mbox{ind\,} \del$.  We clarify
this in  Section~9 by showing how the linearization  $f\mapsto
L_f$ defines a map from $\bM$   to a space ${\cal F}$ of {\em
real} Fredholm operators.  There are natural classes $\kappa_i\in
H^*({\cal F})$, first defined by Koschorke,   that give obstructions to the bundle $\mbox{ind}_\R\,
L$ being an actual vector bundle rather than a virtual bundle. We  prove that
all Koschorke classes vanish for  the family of maps
defined by the moduli space (\ref{M2}), so $\Ob=\mbox{ind}_\R\,
L$  is  an actual bundle over the moduli space.   This gives a homotopy-theoretic characterization of the
obstruction bundle.

 Zinger \cite{Z} also showed (\ref{lp3}) following the approach of Fukaya-Ono \cite{FO} and Li-Tian \cite{LT1}.  He also describes some interesting generalizations of our main theorem.

\vspace{1cm}




\setcounter{equation}{0}
\section{$J_\a$-holomorphic maps with stabilized domains}
\label{section1}
\medskip

We begin with a review of the setup for $J_\alpha$-holomorphic maps
(for details see \cite{L} and \cite{LP}).  Fix a K\"ahler surface
$X$ with complex structure $J$ and geometric genus $p_g>0$.  Then
the  real vector space
\begin{equation*}
\H =\ \mbox{Re}\big(\,H^{2,0}\oplus  H^{0,2}\,\big)
\end{equation*}
has  dimension $2p_g>0$. Using  the K\"{a}hler metric compatible
with $J$,  each $\a\in \H$ defines an endomorphism $K_{\alpha}$ of
$TX$ by the equation
\begin{equation}
\label{defK} \langle u,K_{\alpha}v\rangle=\alpha(u,v).
\end{equation}
Each $K_\a$ is skew-adjoint, anti-commutes with $J$, and satisfies $
K_\a^2\ =\   - |\a|^2 Id.$   It follows that $J K_\a$ is
skew-adjoint and $Id+JK_{\alpha}$ is invertible. Thus there is a
family of almost complex structures
\begin{equation}
\label{jalpha} J_{\alpha} = (Id + JK_{\alpha})^{-1}J\,(Id +
JK_{\alpha})
\end{equation}
on $X$ parameterized by $\a\in \H$.  Note that while $\a$ is
holomorphic, the corresponding almost complex structure $J_\a$ need
not be integrable, and indeed, usually isn't integrable.

 For  each $\a\in \H$,  we can consider  the set of
maps $f:C\to X$ from a connected complex curve with complex
structure $j$ into $X$ that satisfy the $J_\alpha$-holomorphic map
equation
\begin{equation}
\del_{J_\a}f=0 \label{junhoeq2}
\end{equation}
where $\del_{J_\alpha}f=\frac12(df+J_\alpha dfj)$.  This is
equivalent to the ``perturbed $J$-holomorphic map equation''
\begin{equation}
\overline{\partial}_Jf- \nu_\alpha =0 \qquad\mbox{where}\qquad
\nu_\alpha = K_\alpha(\partial_J f)j. \label{junhoeq}
\end{equation}
In this paper, we fix  $J$ and work with  the $J_\alpha$-holomorphic
map equation in the form (\ref{junhoeq}) rather than
(\ref{junhoeq2}).

Because  any holomorphic map into a K\"{a}hler surface represents a
$(1,1)$ class,  we can restrict attention to maps representing (1,1)
classes (the GW invariants vanish for all other classes).  In this
context, the first author observed that the following remarkable
fact.

\begin{lemma}[Image Localization]
\label{ILLemma}
 If $f:C\to X$ is  a  $J_\a$-holomorphic map  that represents a (1,1) class, then the image of $f$   lies in one connected component $D_k$ of the zero set of $\alpha$.
\end{lemma}

This  fact leads to the general formula  (\ref{key}) expressing the
GW invariant as a sum  of  local invariants  associated with the
components of the divisor of  $\alpha$.  When such a divisor $D$ is
smooth with multiplicity one,  the square of its holomorphic normal
bundle $N$ is the canonical bundle $K_D$, so $(N, D)$ is  a spin
curve.  In this case, it was shown in  \cite{LP} that the local GW
invariants depend only on the genus and parity of the spin curve
$(D,N)$ (the parity is $h^0(N)$ mod 2).

\medskip

Now consider   a genus $h$ spin curve $(D,N)$, that is, a curve $D$
with genus $h$ and a holomorphic line bundle $N$ on $D$ with
$N^2=K_D$.   The total space $N_D$  of $N$ has a complex structure
$J$ that makes  the projection $\pi:N_D\to D$ holomorphic.  We then
have an exact sequence
\begin{equation}
\label{decomp} 0 \to \pi^*N \to TN_D \to \pi^*TD \to 0
\end{equation}
 and hence the canonical bundle of $N_D$ is $\pi^*(N^*\otimes K_D)=\pi^*N$.  The tautological section of $\pi^*N$ is thus
 a holomorphic 2-form $\a$ on $N$ that vanishes transversally along the zero section $D\subset N_D$.
 This 2-form $\a$ induces an almost complex structure $J_\a$ on $N_D$ as in (\ref{jalpha}).
 The local GW invariants of $D$ are then associated with the space
\begin{equation}\label{msp-1}
\CM^{J_\alpha}_{g,n}(N_D,d)
\end{equation}
 of all stable $J_\alpha$-holomorphic maps $f:C\to N_D$ whose domain has genus $g$ and $n$  marked points
 and whose image represents $d[D]\in H_2(N_D)$.  By Proposition~\ref{ILLemma} this is the same as the space
 $\CM_{g,n}(D,d)$ of stable  $J$-holomorphic maps into the zero section $D$  with degree $d$.
 However, the linearizations of these equations differ in a way that will be crucial in later sections.

 \medskip

\begin{lemma}\label{lemma-1}
 If  $D$ has genus  $h \geq 1$  and $2g+n \geq 3$, then the domain of every map in $\CM_{g,n}(D,d)$ is a stable curve. 
 \end{lemma}

\pf When $D$ has genus $h \geq 1$ all rational components of the
domain are mapped to points, so are stable curves (by  the
definition of stable map).  Similarly, all components with genus one
are mapped to points or have a node (so are stable), unless the
domain is smooth and the map is etale, in which case the domain has
at least one marked point since $2g+n \geq 3$.  \qed

\medskip

  In particular, because all domains are stable,  the relative cotangent bundles over $\CM_{g,n}(D,d)$ are
pull-backs of the relative cotangent bundles over $\CM_{g,n}$ by the stabilization map.
The descendent classes are thus pull-backs of cohomology classes via the map (\ref{1.evmap}).

\bigskip

One usually takes $\bM_{g,n}$ to be a Deligne-Mumford space.
However, it is more convenient to take it to be
the moduli space described by  Abramovich, Corti and Vistoli \cite{ACV}, building on the key work of Looijenga \cite{Lo}.
This is a finite branched cover of the compactified Deligne-Mumford space that  is a fine moduli space
and a smooth projective variety (Theorem~7.6.4 of \cite{ACV});  it solves the moduli problem for families of  ``$G$-twisted curves'' that carry a principle $G$-bundle for a certain finite group $G$ together with certain additional structure at their nodes and marked points (see \cite{ACV} for details).
 As described in Section~1 of \cite{IP1}, this  modification has no substantial effect:
one recovers the standard GW invariants by dividing by the degree of the cover.
Accordingly, we will use the standard notation $\bM_{g,n}$ for compactified Deligne-Mumford space and
leave the presence of twisted structures implicit.
In this context, the  space $\bM_{g,n}$ of $G$-twisted curves and the total
space of its universal curve
\begin{equation}
\label{Prymfamily}
\begin{CD}
\bU_{g,n}\\
@VV{\pi}V\\
\bM_{g,n}
 \end{CD}
\end{equation}
are manifolds with Riemannian metrics.

\medskip

In the paper we will be working will the moduli  spaces of solutions
of the perturbed $J$-holomorphic map equation
\begin{equation}
\label{3.Jnu} \del_{J_\alpha}f=\nu_f
\end{equation}
for various types of perturbation terms $\nu$.  We will use
perturbations of the type introduced by Ruan and Tian in \cite{RT2},
which  can be described as follows.   Fix an
almost complex manifold $(X,J)$ and
let $\bU$ be
the universal curve (\ref{Prymfamily}).  Because  $G$-twisted curves
have no non-trivial automorphisms, the domain of a map $f:C\to X$ is
uniquely identified (as a $G$-twisted curve) with a fiber of $\bU$, and the
graph of $f$ is a map $F:C\to \bU\times X$.  Consider
the bundle $T^*\bU\boxtimes TX=\Hom(\pi_1^*T\bU,\pi_2^*TX)$ over
$\bU\times X$. A {\em Ruan-Tian perturbation}   is an element $\nu$
of the space
\begin{equation}
\label{defP} {\cal P} = \Omega^{0,1}(T^*\bU\boxtimes TX)
\end{equation}
of $(0, 1)$ sections, that is, sections $\nu$ satisfying  $\nu\circ
j_{u}=-J\circ \nu$ where $j_{u}$ is the complex structure on $\bU$.  Restricting such a $\nu$ to the graph
of $f$ gives a form $\nu_f\in \Omega^{(0,1)}(f^*TX)$, defined by
\begin{equation}\label{restriction}
v_f(x)(u) = v(x,f(x))(u)\ \ \ \ \ \ \ \ \ \ \forall\, x\in C,\ \
u\in T_xC, 
\end{equation}
that can be used in (\ref{3.Jnu}).   Alternatively, one can observe
that $\nu$ defines  an almost complex structure $J_{\nu}$ on
$\bU\times X$ by $J_{\nu}(v,w)=(j_{\bU}(v),  J(w)-2\nu( v))$ and
that
  $f:C\to X$ is $(J,\nu)$-holomorphic  if and only if its graph
$F:C\to \bU \times X$ is $J_\nu$-holomorphic (see \cite{PW} and
\cite{RT1}).

We  will routinely use the phrase ``for generic $\nu$'' to mean ``for $\nu$ in a Baire set in the space ${\cal P}$''.

\vspace{1cm}

\setcounter{equation}{0}
\section{Convergence   of maps with stable domains}
\label{section2}
\medskip

The proof of the main theorem requires a careful definition of the space of maps as a topological space.  The appropriate topology is {\em not}  the one defined by Gromov compactness.  Instead, we will use the stronger ``$\la_p$-topology''  defined carefully-chosen weighted Sobolev norms.
In this section, we set out the definitions and then prove that Gromov compactness holds in the $\la_p$-topology.

Let $\pi:\CU\to \bM_{g,n}$ be the universal $G$-twisted curve (\ref{Prymfamily}).
 There is  an $\delta_0>0$
such that in each  fiber $C_z=\pi^{-1}(z)$ the nodes are separated
from each other and from the marked points by a distance of at least
$4\delta_0$.  We will scale the metric so that $\delta_0=1$ and
choose a function
 \begin{equation}
 \label{defrho}
\rho: \bU\to \R
 \end{equation}
that is equal to the distance to the nodal variety on the set
$\{\rho<1\}$ and satisfies $0\leq\rho\leq 2$ everywhere.  For  each
$\delta<1$ we will also write
\begin{equation}
\label{defCd} C_z(\delta)\ =\ C_z\cap \{\rho\geq  \delta\}
\qquad\mbox{and}\qquad B_z(\delta)\ =\ C_z\cap \{\rho < \delta\}.
\end{equation}
For small $\de$,  $B(\delta)$ is a union of components, each being either a union of two disks with their center points identified (a node) or a thin annular  neck (a near node).  Because the universal $G$-twisted curve (\ref{Prymfamily}) is a fine moduli space,  the fibers have no non-trivial automorphisms and:
\begin{enumerate}
\item[(i)]  Around each smooth fiber $C_z$ there is a neighborhood $V_z$ of $z$ and   local trivialization
 \begin{equation}
 \label{firsttrivialization}
 \phi: C_z\times V_z \to U_z.
 \end{equation}
 \item[(ii)]  Around each fiber $C_z$ with $B_z(\d)\not= \emptyset$  there is a  similar  local trivialization
 \begin{equation}
 \label{secondtrivialization}
 \phi: C_z(\d)\times V_z \to U_z
 \end{equation}
where   $U_z=\pi^{-1}(V_z)\cap \{\rho > \delta\}.$
\end{enumerate}

\medskip

Now fix  an isometric embedding of $X$ into $\R^N$ for some $N$.    The space of maps is defined using the following   weighted Sobolev norms.

\begin{defn}
\label{defn2.1}
For each $p$, $2\leq p < \frac52$, set $\la=\la_p=\frac{p-2}{6}$.  For
each curve $C$, let $\Map_\la(C,X)$ be the completion of the set of maps $f:C\to X$ in the
weighted norm
  \begin{equation}
 \label{alphanorm}
 \|f\|^p_{1, p, \la}\ =\ \int_C  \rho^{-\la}\, |df|^p\ +\ \rho^{\la-2} \, |f|^p\ dvol_C
 \end{equation}
where $|f|$ is defined by the embedding $X\subset \R^N$.  
 \end{defn}

 If we replace $C$ by $C(\delta)$ in the Definition~\ref{defn2.1} then the resulting
 norms are uniformly equivalent,  for each $\delta$, to the usual  $W^{1,p}$ norm on the fibers of  the local trivializations
 (\ref{firsttrivialization}) and   (\ref{secondtrivialization}). We also define the {\em $p$-energy} of a map $f:C\to X$ to be
 \begin{equation}
 \label{pEnergy}
E_p(f)\ =\   \left(\int_C \rho^{-\la}\, |df|^p\right)^{\frac{2}{p}}.
 \end{equation}
This is the usual energy of $f$ when $p=2$. We will consider only
the subset of $\Map(X)$ whose  $p$-energy is  below a fixed level.

\begin{defn}[$\la_p$-topology]
\label{defMapE} For each number $E$ and $p>2$, set
\begin{equation}
\label{MapE} \Map_{g,n}^{E}(X)\ =\ \left\{ (z,f)\,|\, z\in
\CM_{g,n},\ f\in \Map_\la(C_z,X)\ \mbox{and
}\ E_p(f) < E\, \right\}.
\end{equation}
Give this space the $\la_p$-topology:    a sequence $(C_n, f_n)$
converges to $(C,f)$  if (a)\, $C_n\to C$ in
$\bM_{g,n}$,  (b)\, $E_p(f_n)\to E_p(f)$,  and (c)\, $f_n \to f$ in
the norm (\ref{alphanorm}) on $C(\delta)$ for every $\delta>0$.
\end{defn}

The convergence in (c) is defined after identifying the domains
using  the trivialization (\ref{secondtrivialization}).

\medskip

When working with the $\la$-norm on nodal curves $C$, it is helpful to use the
definitions  (\ref{defCd}) to decompose $C$  into the part
$C(\delta)$ away from the nodes, and the neighborhoods $B_i(\delta)$
of the nodes $n_i$ and to parameterize each $B_i(\delta)$ by a
cylinder as follows.  Near $n_i$, choose holomorphic coordinates
$\{x^i\}$ in the universal curve such that $B_i(\delta)$ is $\{
(x^1, \dots x^n)\, |\ x^1x^2=\mu \}$.  Set $L(\mu)=|\,\ln
\delta-\frac12 \ln |\mu|\,|$ and map the cylinder $T(\delta)=[-
L(\mu), L(\mu)] \times S^1$  to $B_i(\delta)$ by
  \begin{equation}
 \label{defphi}
\phi_\mu:  (t,\theta)\mapsto \left(x,  \frac{\mu}{x}\right) \quad
\mbox{where}\quad  x=\sqrt{|\mu|} \,e^{t+i\theta}.
 \end{equation}
Then the pullback of $\rho^2$ is $2|\mu|\cosh (2t)$ and the metric
$g$ on $C$ is conformally related to the metric $\hat{g}$ on the
cylinder by $\phi^*g=\rho^2\, \hat{g}$.  Consequently, after
identifying $f$ and $\rho$ with their pullbacks,
  \begin{equation}
 \label{alphanorm2}
 \int_{B_z(\delta)}  \rho^{-\la}\, |df|^p\ +\ \rho^{\la-2} \, |f|^p\ dvol_g\ =\  \int_{T(\delta)}  \rho^{ -7\la}\, |df|^p\ +\  \rho^{\la} \, |f|^p\ dvol_{\hat{g}}.
 \end{equation}
 (The exponent $7\la$ arises because $|df|^p\, dvol_g$ pulls back to $ \rho^{2-p}\, |df|^p\, dvol_{\hat{g}}$, and $2-p-\la=-7\la$.)
For $\mu=0$, we can replace  $T(\de)$ by two cylinders, both isometric to  $[0,\infty)\times S^1$,  and replace (\ref{defphi})  by the maps $(t,\theta) \mapsto (e^{t+i\theta}, 0)$ and $(s,\theta') \mapsto (0, e^{s+i\theta'})$. Then (\ref{alphanorm2}) holds for all $\mu$.   Whenever we mention  (\ref{defphi})   and (\ref{alphanorm2}) it should be understood that the $\mu=0$ case is defined in this way.  Because $\rho$ is essentially exponential in $t$,   the integrals of its powers in the cylindrical metric  satisfy 
\bear
\int_{\rho\leq\delta} \rho^\gamma\ dt\,d\theta\ \leq\ c_\gamma\,
\delta^\gamma \quad \mbox{and}\ \quad \int_{\rho\geq\delta}
\rho^{-\gamma}\ dt\,d\theta\ \   \leq\ c_\gamma\, \delta^{-\gamma}
\hskip1.6cm \mbox{for}\    \gamma> 0.
 \label{integralofrho}
 \eear

 \medskip

 The following two lemmas give properties of maps in $\Map_{g,n}^{E}(X)$.  The first shows that the images of balls around nodes are uniformly small, and the second shows that when two maps are close in the $\la_p$-topology then the Hausdorff distance $\dist_{\H}$ between their images is small.

\begin{lemma}
\label{diamlemma}
Fix $g,n, E$ and $p>2$.  Then  there are constants $\delta_0>0$   and  $c$ such that for every $(C, f)\in \Map_{g,n}^E(X)$ and every node $n_i$ of $C$ we have, for all $\delta<\delta_0$,
$$
\mbox{\rm diam}(f(B_i(\delta)))\ \leq\ c\,  \delta^{\frac{7\la}{p}}.
$$
\end{lemma}

\pf   Fix  $\delta_0$ small enough that the maps (\ref{defphi}) are defined on the $2\delta_0$-neighborhoods of the nodes. For each node $n_i$ set $\Omega=\phi^{-1}(B_i(\delta))$.   Fix some $q$ with $2<q<p$.  On the cylinder $T$,  the oscillation of $f$ is bounded by the $L^{q}$ norm of $df$:
$$
\mbox{\rm diam}\, f(B_i(\delta))\ =\
\mbox{\rm diam}\, f(\Omega)\ = \
\underset{\Omega}{\mbox{\rm osc}}\, f\ \leq\
c_2 \left(\int_{\Omega'} |df|^{q}\right)^{\frac{1}{q}}.
$$
where $\Omega'=\phi^{-1}(B_i(2\delta))$.   The proof is completed by inserting
 $1=\rho^{-7\la q/p}\rho^{7\la q/p}$ into the integrand, applying
Holder's inequality  and then using (\ref{alphanorm2}) and (\ref{integralofrho}):
\begin{equation}
\label{2.holder}\notag
\mbox{\rm diam}\, f(B_i(\delta))\ \leq\
c_3\,   \delta^{\frac{7\la}{p}}\, \left(\int_{\Omega'} \rho^{-7\la}|df|^{p}\right)^{\frac{1}{p}}  \ \leq\
c_4\,  \delta^{\frac{7\la}{p}}\, \sqrt{E_p(f)}\ \leq\
c_5 \  \delta^{\frac{7\la}{p}}.
\end{equation}
\qed

\begin{lemma}
\label{Hdistlemma}
Fix $g,n, E$ and $p>2$.  Given $(C,f)\in \Map_{g,n}^E(X)$, $\ep>0$ and  a sufficiently small $\delta>0$, there are
neighborhoods $\N_f(\ep)$  and $\N_f(\ep, \delta)$ of $(C,f)$ in the $\la_p$-topology such that
\begin{eqnarray*}
\dist_\H(f(C), g(C')) \leq \ep  \qquad \mbox{for all $(C', g)\in\N_f(\ep)$}
\end{eqnarray*}
and
$$
 \|f-g\|_{\infty; C(\delta)}\leq \ep \qquad \mbox{for all $(C', g)\in\N_f(\ep, \delta)$.}
$$
\end{lemma}

\pf  By Lemma~\ref{diamlemma} we can choose $\delta_1>0$ so that $\mbox{\rm diam}(f(B_i(\delta_1))) \leq \ep/3$ for all $f\in \Map_{g,n}^E(X)$.    Fix $\delta$ with $2\delta\leq \delta_1$.  We can then choose a neighborhood $\N_f$ of $(C,f)$ small enough that the domains of all maps $g:C'\to X$ in $\N_f$ lie in the uniform local trivialization (\ref{secondtrivialization}) around $C$ for this  $\delta$.  Decompose the domain $C$ of $f$ into $C(\delta)$ and the union of  neighborhoods $B_i(\delta)$ of its nodes.  
Then for each $(C', g)\in\N_f$ we have  $\delta\leq \rho\leq 2$ on $C'(\delta)$, so the norm (\ref{alphanorm}) is uniformly equivalent to the (unweighted) $W^{1,p}$ norm on $C'(\delta)$.  Furthermore, the fibers of the  local trivialization (\ref{secondtrivialization}) have uniform geometry, so there is a uniform constant (depending on $\delta$) for the usual Sobolev embedding $C^0\subset W^{1,p}$.  Consequently,
$$
 \|f-g\|_{\infty; C(\delta)} \ \leq\   c(\delta)\,\|f-g\|_{1,p, \la}.
$$
We can then choose $\N_f(\ep, \delta)$ to make the right-hand side less that $\ep/3$ for all $g\in \N_f(ep, \delta)$.  In particular,  setting $\N_f(\ep)=\N_f(\ep, \delta_1/2)$, each $g\in \N_f(\ep)$ satisfies
\begin{eqnarray*}
\dist_\H(f(C), g(C'))& \leq &   \|f-g\|_{\infty; C(\delta)} \,+\, \sup_i \mbox{\rm diam}\, f(B_i(\delta))  \,+\, \sup_i \mbox{\rm diam}\, g(B_i(\delta))\\
& \leq & \frac{\ep}{3} +  \frac{\ep}{3} +  \frac{\ep}{3}.
\end{eqnarray*}
\qed

\medskip

We next prove  an enhanced version of the Gromov Compactness Theorem. It assumes that all maps  have stable domains (cf. Lemma~\ref{lemma-1}), and proves convergence in the $\la_p$-topology, which is stronger than the convergence is the standard formulations of Gromov Compactness.

\begin{theorem}[Compactness]
\label{compactnesstheorem}
Suppose that  all maps in $\bM_{g,n}^{J,\nu}(X,A)$ have stable domains.  Then
there is an $E=E(p,g,n,A)$ such that  $\bM_{g,n}^{J,\nu}(X,A)$ is a compact subset of $\Map^E(X)$ whenever $\sup \nu$ is
small.
\end{theorem}

\pf  Given a sequence of    maps $f_k:C_k\to X$  in $\bM_{g,n}^{(J,\nu)}(X,A)$ we can consider their graphs
 $F_k:C_k \to \bU \times X$ as described after equation~(\ref{restriction}).  Then $\{F_k\}$ is a sequence of $J_\nu$-holomorphic maps with uniformly bounded energy.  Moreover, as explained in Section~1,  each $C_k$ is a stable curve with no non-trivial automorphisms.  The first factor of each $F_k$ is therefore a diffeomorphism onto a fiber of the universal curve.
The  Gromov Compactness Theorem \cite{IS}, applied to $\{F_k\}$  implies that there is a subsequence such
 that (i) the domains $C_k$ converge in $\bU$ to a limit $C_0$, 
(ii) the maps converge in Hausdorff distance and in
$W^{1,p}$ on compact sets in the complement of the nodes  of $C_0$, and (iii) the energy densities $|dF_k|\,dvol$  converge as measures.

Next  fix $\delta$ small enough that the maps (\ref{defphi}) are defined on the $2\delta$-neighborhoods of the nodes.  Because the  energy densities converge as measures we can, given any $\ep_0>0$,  also assume the energy in the neighborhoods $B_z(\delta)$ of each node $z$ in each $C_k$ is at most $\ep_0$.  But energy is conformally invariant so, as in (\ref{alphanorm2}), the pullback maps $f_k:T(\delta) \to X$ satisfy the $J$-holomorphic map equation and
 \begin{equation}
 \label{ep0bound}
\int_{T(\delta)} |df_k|^2\ =\ \int_{B_z(\delta)} |df_k|^2\ <\ \ep_0.
 \end{equation}
Elliptic theory then gives a pointwise bound on $|df_k|$:   by Lemma~5.1 of \cite{IP2} there are universal constants $c$ and  $\ep_0$ so that whenever (\ref{ep0bound}) holds we have $|df_k|<c\rho^{1/3}$ on $T(\delta)$ where $c$ and the maximum number of nodes  depend only on $g,n,A$ and $\sup |\nu|$.  Integrating by (\ref{alphanorm2}) and  (\ref{integralofrho}) shows that, for $2<p<5/2$,
 \begin{equation}
 \label{aboundonM}
 \int_{B(\delta)}  \rho^{-\la}\, |df_k|^p\ =\ \int_{T(\delta)}  \rho^{ -7\la}\, |df_k|^p\ \leq\  c \int_{T(\delta)}  \rho^{-7\la +\frac{p}{3}}\ \leq\ c\,\delta^{\frac{1}{4}}.
 \end{equation}
 Convergence in the $\la$-norm follows:  given $\ep>0$,  choose $\delta$ small enough that the above bound is less than $(\ep/2)^p$ and then choose $K$ large enough that, for all $k\geq K$, we have $\|f_k-f_0\|_p < \ep/2$  on the complement $C(\delta)$ of $B(\delta)$.  Finally, note that (\ref{aboundonM})    bounds $E_p(f)$   for all $(J,\nu)$-holomorphic maps $f$.
 \qed

\medskip

\vspace{1cm}

\setcounter{equation}{0}
\section{The virtual fundamental class}
\label{section3}
\medskip

Gromov-Witten invariants of a closed symplectic manifold $X$ have been defined in several different ways.
Some of these have been proved to be equal, and it is folklore that all define equivalent invariants.
Each involves compact  moduli spaces $\CM_{g,n}(X,A)$ of (formal) dimension
$r=c_1(X)(A) + (\dim X-3)(1-g) + n$ and the associated  map
\begin{equation}
\label{9.evmap}
\hev=st\times \ev:\CM_{g,n}(X,A) \to \bM_{g,n}\times X^n
\end{equation}
defined as in (\ref{1.evmap}).   For our purposes,  three
different descriptions are relevant.  One is algebraic and the other
two are analytic:

\begin{enumerate}
\item When $X$ is K\"{a}hler, the space of stable maps is a projective variety
$\bM^{\mbox{\rm \scriptsize stable}}_{g,n}(X,A)$ and algebraic
geometers  (see \cite{LT2} and \cite{BF}) define the virtual
fundamental class $[\bM_{g,n}(X,A)]^{virt}$ as an element of its
Chow cohomology:
\begin{equation}
\label{Chowclass} [\bM_{g,n}(X,A)]^{virt}\,\in\, A_r(\bM^{\mbox{\rm
\scriptsize stable}}_{g,n}(X,A)).
\end{equation}

\item When $X$ is  symplectic and semipositive,  Ruan-Tian \cite{RT2} showed that for generic $(J,\nu)$
the space $\bM^{(J,\nu)}_{g,n}(X,A)$ of all $(J,\nu)$-holomorphic maps is the union  of
(i) a  $2r$-dimensional  orbifold $\M$
 consisting of maps from smooth domains  and (ii) a stratified ``boundary''  whose image under the map
 (\ref{9.evmap})  lies in a set of dimension    at most $2r-2$.
Consequently, the image of the moduli space represents a rational
homology class
\begin{equation}
\label{GW-class} GW_{g,n}(X, A)\,\in\,H_{2r}(\CM_{g,n}\times
X^n;\Q).
\end{equation}
This Ruan-Tian GW class  is independent of  the generic $(J,\nu)$ and is a symplectic invariant of $X$.
(``Semipositive'' is a technical condition that is true whenever $\dim X\leq 6$.)  \\

\item When $X$ is  symplectic,  the construction  of Li-Tian \cite{LT1}   defines a virtual fundamental class
\begin{equation}\label{LTVFC}
[\CM_{g,n}(X,A)]^{\vir}\,\in\,H_{2r}(\Map_{g,n}(X,A);\Q).
\end{equation}
in the homology of the infinite-dimensional  space of maps
$\Map_{g,n}(X,A)$.  Because (\ref{9.evmap}) extends to $\Map_{g,n}(X,A)$
one can then evaluate the class (\ref{LTVFC}) on
classes in $H^*(\CM_{g,n}\times X^n)$.

\end{enumerate}

Variations on construction 3 have been done by Fukaya-Ono \cite{FO},
Ruan \cite{R} and Siebert \cite{S}.

\medskip
When $X$ is  K\"{a}hler,  Li and Tian proved in \cite{LT3}  that the
virtual class (\ref{LTVFC}) is the homology class underlying the
Chow class (\ref{Chowclass}) under the inclusion of the space of
stable maps into  $\Map_{g,n}(X,A)$.  When $X$ is semipositive,  one can show    that  the
pushforward of the virtual class (\ref{LTVFC})  by $\hev_*$
is the  Ruan-Tian class (\ref{GWclass});  we explicitly prove this for compact curves in Remark~\ref{RT=LT}.  Thus, when $X$
is a smooth compact curve $D$,  all three definitions apply and, after  pushing forwards
by $\hev_*$, define the same element in the homology of
$\CM_{g,n}\times X^n$.

\bigskip\medskip

Later, to prove the Main Theorem stated in the introduction, we will use a variation of  Definition (\ref{GW-class}).
We only consider the case relevant to the present paper: when the domains of all $(J,\nu)$-holomorphic maps
have no non-trivial automorphisms (cf. Section~1).
In this case, for generic $(J,\nu)$,  each stratum in the moduli space
\begin{equation}
\label{3.ms}
\CM\,=\,\bM^{(J,\nu)}_{g,n}(X,A)
\end{equation}
(including the ``top stratum''  $\M$)  is a smooth oriented manifold,
 the boundary $\bd\CM=\CM\ssetminus \M$
is compact,  and $\hev$ restricts to a smooth map to the compact oriented manifold $\CM_{g,n}\times X^n$
on each stratum.   The Ruan-Tian invariant
 can be defined, using arguments of Kronheimer and Mrowka,  by deleting a neighborhood of the boundary and considering the resulting relative homology class, as follows.

First, a simple compactness and transversality  argument (see the first paragraph of  the proof of Proposition~4.2 in \cite{KM}) shows that
there is an open neighborhood $U$ of
the image $\hat{ev}(\bd\CM)$ satisfying\,:
\begin{itemize}
\item[(a)]
$\overline{U}$ is a smooth manifold with boundary and
$\hat{ev}_{|\M}$  is transverse to $\bd U$, and
\item[(b)]
there is a (finite) basis for $H_{\ell-2r}(\CM_{g,n}\times X^n;\Q)$
represented by cycles disjoint from $\overline{U}$
\end{itemize}
where $\ell$ is the dimension of $\CM_{g,n}\times X^n$.

\medskip

\begin{defn}
\label{defn3.1}
Given a moduli space (\ref{3.ms}), choose a set $U$ as above and let $\bM_U\subset\,\M$ be the closed set
\begin{equation}\label{rel-moduli}
\bM_U\,=\,  \CM\ \,\cap\, \hat{ev}^{-1}(\CM_{g,n}\times X^n\ssetminus U)\,\subset\,\M.
\end{equation}
It follows from (a) above that $\bM_U$ is a compact oriented
manifold with boundary, and thus carries a fundamental class in
relative homology:  $[\CM_U]\in H_{2r}(\M_U,\bd\M_U; \Q)$.
\end{defn}

On the other hand, the inclusion of the pairs $k: (\CM_{g,n}\times X^n,\emptyset) \to (\CM_{g,n}\times X^n,\overline{U})$
induces a map in homology
$$
k_*:\,H_{2r}(\CM_{g,n}\times X^n;\Q) \to H_{2r}(\CM_{g,n}\times X^n,\overline{U};\Q).
$$
Proposition 4.2 of \cite{KM} then implies that there is a unique homology class $[\hev(\CM)]$ with
\begin{equation}\label{unique-class}
k_*[\hev(\CM)]\,=\,\hev_*[\CM_U].
\end{equation}
In fact, they showed that there is a rational singular smooth cycle $B$ that represents $[\hev(\CM)]$ and
agrees with $\hev(\CM_U)$ outside of $U$. The uniqueness then follows from property (b) above because
the inclusion
$$
H_{2r}(\overline{U};\Q)\to H_{2r}(\CM_{g,n}\times X^n;\Q)$$
is trivial and hence $k_*$ is injective.
The class $[\hev(\CM)]=[B]$ is the GW class\,:

\begin{lemma}
\label{KMGW}
With $U$ and $\bM_U$ as above, the Ruan-Tian GW class is the unique rational homology class satisfying (\ref{unique-class}), namely
\begin{equation}\label{GWclass}
GW_{g,n}(X,A)\, =\,[\hev(\CM)]\,\in\,H_{2r}(\CM_{g,n}\times X^n;\Q).
\end{equation}
In particular, $[\hev(\CM)]$ is independent of the choice of $U$ and the generic $(J,\nu)$.
\end{lemma}

\pf  For each $\gamma\in H^{2r}(\CM_{g,n}\times X^n;\Q)$, Ruan and Tian choose a representative $\Gamma$ of
the Poincar\'e dual of $\gamma$ with $\Gamma$ transversal to $\hev(\bM)$ and disjoint from $\hev(\bd\CM)$.
They defined $GW_{g,n}(X,A)(\gamma)$ as $\hev(\M)\cap \Gamma$ and showed these numbers are independent of
the various choices made.  But by property (b) above, we can assume that $\Gamma$ is a linear combination of
cycles that do not intersect $\overline{U}$.  Then the geometric intersection $\hev(\M)\cap \Gamma$ is the same as $B\cap \Gamma$,
which represents $[B](\gamma)=[\hev(\CM)](\gamma)$.
\qed

\medskip

In our case,  $X=N_D$ is the total space of the spin curve $(D,N)$.
By the image localization Lemma~\ref{ILLemma},
for small generic $\nu$ there is a $\epsilon$-neighborhood $N_D(\ep)$ of the zero section $D$ in $N_D$ such that
every $(J_\a,\nu)$-holomorphic map has its image in that neighborhood. Although $N_D(\ep)$  is not compact, it is a manifold with boundary and  Lemma~\ref{KMGW} still applies (again by Proposition~4.2 of \cite{KM}).
Thus the image of the moduli space represents a rational homology class
\begin{equation}\label{LGWclass}
GW_{g,n}^{loc}(N_D,d)\,:=\,\big[\hev\big(\,\bM^{(J_\a,\nu)}_{g,n}(N_D,d[D])\,\big)\big]\,\in\,
H_{2\b+2n}(\CM_{g,n}\times N_D^n;\Q).
\end{equation}
Here the number $\b$, which  will occur frequently in our dimension counts, is given in terms of the degree $d$ of the map and the genus $h$ of $D$ by
\begin{equation}\label{defbeta}
\b\,=\,d(1-h)+g-1.
\end{equation}
The class (\ref{LGWclass}) is the local GW class of the spin curve $(D,N)$.
As shown in \cite{LP}, for given $g,n,d$ and $h$, it
depends only on the parity $h^0(N)$.

 \vspace{1cm}

\setcounter{equation}{0}
\section{The linearization operator}
\label{section4}
\medskip

We now return to the specific situation in which the local GW
invariants are defined.  Thus, as described after Lemma~1.1, we
consider $J_\a$-holomorphic maps into the total space $N_D$ of a
spin bundle $N\to D$ over a curve $D$.    This section shows how the
special form of the $J_\a$-holomorphic map equation implies
vanishing theorems for the linearized operator.

\bigskip
A neighborhood $U$ of the zero section  $D\subset N_D$ can be
identified with   a neighborhood of the zero section of the
projectivization $\P(N\oplus \O_D)$ by a fiber- preserving
biholomorphism.  Since $\P(N\oplus \O_D)$ is
  K\"{a}hler we can use this identification to define a K\"{a}hler
structure  on $U$. With this K\"{a}hler structure, the exact
sequence (\ref{decomp}) of the underlying complex vector bundles
splits as
\begin{equation}\label{split}
TN_D\ =\ \pi^*TD\oplus \pi^*N.
\end{equation}
We will often write  this simply as  $TN_D=TD\oplus N$.  In this and
the following sections $\nabla$  will denote the normal component of
the  Levi-Civita connection of the K\"{a}hler metric  on $N_D$, that
is, $\piN\circ \nabla^{LC}\circ \piN$ where $\pi_N$ is the
orthogonal projection onto the normal part of (\ref{split}).

Let $f:C\to N_D$ be a map  from a smooth domain $C$. The
linearization  of the $J_\alpha$-holomorphic map equation is the operator
$$
D_f:\Gamma (f^*TN_D)\oplus H^{0,1}(C,TC) \to
\Gamma(\Lambda^{0,1}(f^*TN_D))
$$
 given by $D_f(\xi,k)=\widehat{L}_f \xi+J_\alpha df k$;  the operator $\widehat{L}_f$ arises from the variation in the map with the complex structure on $C$ held fixed, and the term  $J_\alpha df k$ arises from the variation of the
 complex structure on $C$.    Under the splitting (\ref{split}), $\widehat{L}_f$ decomposes as
 \begin{equation*}
\widehat{L}_f=\begin{pmatrix}\del^T_f & A \\ B &  L_f\end{pmatrix} :
\Omega^0(f^*TD)\oplus \Omega^0( f^*N) \to
                            \Omega^{0,1}(f^*TD)\oplus \Omega^{0,1}( f^*N)
\end{equation*}
where $A$ and $B$ are bundle maps {\em that vanish when the image of $f$ lies
in the zero section of $N_D$}.
The normal component is given explicitly by
\begin{equation}
\label{Lfformula}
 L_f \xi (v) = \piN\left[\nabla_{v} \xi +
J\nabla_{jv} \xi  - \nabla_\xi K_\a (df(jv)+Jdf(v))-
K_\a(\nabla_{jv} \xi + J\nabla_{v} \xi )\right]
\end{equation}
where $\xi\in \Gamma(f^*N)$ and $v\in TC$.
If $f$ is $J_\a$-holomorphic then this operator reduces to
\begin{equation}
\label{LRformula}
L_f\,=\,\del_f +R_\a
\end{equation}
where $\del_f$, given by the first two terms of (\ref{Lfformula}), is the $\del$-operator of the holomorphic bundle $N\to D$  and $R_\a:f^*N\to T^*C\otimes f^*N$ is the complex anti-linear bundle map
defined by $R_\a(\xi)=-2\nabla_\xi K_\a\circ df\circ j$.

For  a map  $f:C\to N_D$ whose domain $C$ is a connected nodal curve, $L_f$ can be described in terms of the 
normalization $\pi:\tilde{C}\to C$.  The inverse image of each node $n_i\in C$ is  a pair of points $p_i, q_i\in \tilde{C}$.   For each component $\tilde{C}_k$
of the normalization, let $E_{k, f}$ be   the space $\Omega^{0}(\tilde{C}_k,
\pi^*f^*N)$ of smooth sections of $\pi^*f^*N$ on $\tilde{C}_k $, and similarly let
$F_{k,f}=\Omega^{0,1}(\tilde{C}_k,  \pi^*f^*N)$.   Combine these by setting
\vspace{-.2cm}
\begin{equation}
\label{DefEF} 
E_f = \Big\{\xi\in \underset{k}{{\textstyle\bigoplus}} \,E_{k,f}\ | \
\xi(p_i)=\xi(q_i) \ \text{for all $i$}\Big\}\ \ \ \ \ \ \
\mbox{and}\ \ \ \ \ \ \
F_f= \underset{k}{{\textstyle\bigoplus}} F_{k,f}.
\end{equation}
The normal component of the
linearization of the $J_\alpha$-holomorphic map equation is then identified with the 
operator
\begin{equation}
\label{Lbounded} L_f:  E_f  \to F_f
\end{equation}
whose restriction to each component $\tilde{C}_k$ gives  operators $L_{k,f}: E_{k,f} \to F_{k,f}$   as in (\ref{Lfformula}).
Using  Riemann-Roch (cf. Lemma 12.2 of \cite{FO}) one obtains
\begin{equation}
\label{index-nodal} \ind L_f \,=\, -2\beta
\end{equation}
where $\b$ is given by (\ref{defbeta}).

\medskip
In general, as $f$ varies over the space of stable maps, one expects
the dimensions of the kernels of these operators to jump, with
compensating jumps in the dimensions of the cokernels.  But the
following theorem shows that this does not happen for the operators
$L_f$.  This  second remarkable fact about the
$J_\alpha$-holomorphic map equation plays a crucial role in our
analysis.    It implies, as we will show in Sections~5 and 6,  that cokernels of the operators
$L_f $ form vector bundles over $\CM_{g,n}(D,d)$.

\begin{Vthm}
\label{Vtheorem} For each map $f:C\to D$ in $\CM_{g,n}(D,d)$ with
$d\neq 0$
$$
\ker L_f \,=\, 0\ \ \ \ \ \ \
\mbox{and}\ \ \ \ \ \ \
\dim\!\cok L_f\,=\,2\b
$$
\end{Vthm}
\pf  Let $d_k$ be the degree the restriction of $f$ to one component $\tilde{C}_k$ of the normalization of $C$.
If $d_k\ne 0$, $L_{k,f}$ is injective by Proposition 8.6 of \cite{LP}.
On the other hand, when $d_k=0$, the operator $L_{k,f}$ is the $\del$-operator on the trivial bundle
whose kernel is the constant functions. Since any solution of $L_f\xi=0$ restricts to a solution of $L_{k,f}\xi=0$
on  $\tilde{C}_k$, $\xi$ vanishes on each component with $d_k\ne 0$ and is constant on each component with $d_k=0$.
But $\xi$ is continuous at each node, so $\xi\equiv 0$. This shows $L_f$ is injective and hence the dimension of its
cokernel is the negative of its index (\ref{index-nodal}).
\qed

\medskip

For maps $(C,f)$ near a map $(C_0, f_0)$ with nodal domain, it is useful to decompose $E_f$   into a subspace $E^0_f$ and a finite-dimensional complementary subspace $\bar{E}_f$ in such a way that, when $C$ is near $C_0$ and has the same number of nodes,  
\begin{equation}
\label{defE^0}
E_f^0\ =\ \{\xi\in E_f\  |\ \mbox{$\xi(n_i)=0$ at each node $n_i$}\}
\qquad\mbox{\rm and}\qquad
\bar{E}_f \cong \bigoplus_i (f^*N)_{n_i}.
\end{equation}
First, use parallel transport along radial geodesics in $N_D$  to trivialize $N\to N_D$ over a neighborhood $U_i$ of the image $f_0(n_i)$ of each node of $C_0$ (these $U_i$ may overlap).    For small $\delta$, this  trivializes $f^*N$ over the necks $B_i(\delta)$ of all maps $(C,f)$ sufficiently close to $(C_0, f_0)$, and under this trivialization each $v_i\in(f^*N)_{n_i}$ extends to a ``constant'' section $\overline{v_i}$ on $B_i(\de)$.  Next, fix a smooth non-increasing function $\beta(t)$ with $\beta=1$ for $t\leq \frac12$, $\beta=0$ for $t\geq 1$.  By composing with (\ref{defrho}) we obtain, for  each $\de$, a bump function $\beta_\de=\beta(\rho/\de)$ on the universal curve that is supported in a $\de$-neighborhood of the nodal set.  The    subspace
  \begin{equation}
\label{defEbar}
\bar{E}_f   \,=\,   \mbox{span}\big\{\beta_\de \overline{w_i}\ |\, w_i\in(f^*N)_{n_i}\big\}  \subset E_f
\end{equation}
 is a finite-dimensional vector space of smooth sections. Evaluation at the nodes gives the isomorphism in (\ref{defE^0}).  Define a projection $E_f\to \bar{E}_f$, written $\xi\mapsto\bar{\xi}_\de$,   by first expanding $\xi=\sum \varphi_i\overline{w_i}$ for functions $\varphi_i$ on $B(n_i, \de)$ and  then setting
$$
\bar{\xi}_\de=\sum_k  a_i \,\beta_\de\, \overline{w_i}\qquad\mbox{\rm where}\qquad
a_i\, =\, \big(\mbox{vol\,$B(n_i, \delta)$}\big)^{-1}\int_{B(n_i, \de)}\varphi_i
$$
In this context,  each $\xi\in\E_f$ can be uniquely written as the sum $\xi=\bar{\xi}_\de +\zeta_\de$ where
  \begin{equation}
\label{defxibar}
\zeta_\de = \xi-\bar{\xi}_\de.
\end{equation}
vanishes at each node of every nodal curve near $C_0$, so is an element of $E_f^0$.   In the next section we will use this decomposition to complete $E_f$ in a Sobolev norm.

\vspace{1cm}

\setcounter{equation}{0}
\section{Estimates on the linearization}
\label{section5}
\medskip

We now pass to the level of global analysis by completing the spaces $E_f$
and $F_f$  in (\ref{Lbounded}) in appropriate Sobolev norms.  As $f$ varies,
these define vector bundles $\E$ and $\F$ over the space of maps and
the  linearization (\ref{Lbounded}) defines a vector bundle map
$L:\E\to\F$.  However, while $\E$ and $\F$ are topological vector
bundles as in  \cite{LT1}, they are {\it not locally trivial bundles } over maps
with nodal domains.  We will address the issue of local triviality
in the next section. Here, in preparation, we define weighted
Sobolev norms and  establish some elliptic estimates for $L_f$.   Some care  is needed because on weighted Sobolev spaces the Rellich compactness lemma may fail, and   elliptic operators need not be Fredholm  (cf. \cite{lockhart}).

In Section~2, we fixed a Riemannian metric on the universal
$G$-twisted curve and a defining function $\rho$ for the nodal set.  For each $p\geq 2$ we again  set $\la=
\frac{p-2}{6}$. Given a map $f:C\to N_D$ we use $\rho$ to define weighted Sobolev spaces $\E_f$ and $\F_f$ and a decomposition $\E_f\ =\ \E_f^0\oplus \bar{E}_f$:
\begin{itemize}
\item  Let $\F_f$ be the completion of the space of $\Omega^{0,1}(C,f^*N)$  in the weighted $L^p$  norm
\begin{equation}
\label{pnorm} \|\eta\|_{p,\la}^p\ =\  \int_C \rho^{-\la}  |\eta|^p\ dvol_C.
\end{equation}

\item  Similarly, let $E^0_f$ be the space defined in (\ref{defE^0}) and let   $\E_f^0$ be its completion   in the weighted Sobolev norm
\begin{equation}
\label{1pnorm} 
\|\zeta\|_{1, p,\la}^p\ =\  \int_C\rho^{-\la}\  |\nabla
\zeta|^p \ +\ \rho^{-\la-2} |\zeta|^p\  dvol_C.
\end{equation}
This  is well-defined because any smooth $\zeta$ that vanishes at the nodes satisfies $|\zeta|\leq c\rho$ pointwise, and hence the lower-order term in (\ref{1pnorm}) is finite.
\item Let $\bar{E}_f$   be the finite-dimensional space (\ref{defEbar}) for some $\de$ (a specific $\de$ will be chosen in Proposition~\ref{basicellipticlemma} below). The hermitian metric on $N$ induces a metric on $\bar{E}_f$ through the isomorphism (\ref{defE^0}), namely $\|\bar{\xi}\|^2= \sum \left|\xi(n_i)\right|^2$.  

\item  Finally, write $E_f=E_f^0\oplus \bar{E}_f$ as in (\ref{defE^0}) and set
$\E_f\ =\ \E_f^0\oplus \bar{E}_f$. Thus the norm on $\E_f$ is defined by  writing $\xi=\zeta_{\de} +\bar{\xi}_{\de}$ as in (\ref{defxibar}), and setting
\begin{equation}
\label{sumnorm}
\|\xi\|_{1,p,\la}\ =\ \|\zeta_{\de} \|_{1,p,\la} + \|\bar{\xi}_{\de}\|
\end{equation}
\end{itemize}

Notice that in any neck region $B(\delta)$ we can rewrite (\ref{1pnorm}) as an integral over the cylinder $T(\delta)$ using the conformal transformation (\ref{defphi}):
\begin{equation}
\label{1pnormB} 
\|\zeta\|_{1, p,\la; B(\delta)}^p\ =\  \int_{T(\de)}\rho^{-7\la} \left(|\nabla
\zeta|^p \ +\ |\zeta|^p\right)\  dvol_C.
\end{equation}

\subsection{Bounds on $L_f$}

 When $p=2$    the above norms reduce to the usual
$W^{1,2}$ and $L^2$ norms.   When  $p>2$, (\ref{1pnorm})  dominates the
$L^\infty$ norm:

\begin{lemma}
\label{supbound} For each $p>2$ there is a constant $c$, depending
on $p$ but uniform for domains $C$ in the universal curve and the
map $f$,  such that
\begin{equation}
\label{supboundeq}
 \|\xi\|_\infty\leq c\, \|\xi\|_{1,p,\la} \qquad
\mbox{for all } \xi\in \E_f.
\end{equation}
\end{lemma}

\pf   Fix  $\delta_0$ small enough that the maps (\ref{defphi}) are
defined on the $2\delta_0$-neighborhoods of the nodes.   Each fiber
$C_z$ of the universal curve decomposes as $C_z(\delta)\cup
B_z(\delta)$ as in (\ref {defCd}). The  domains $C_z(\delta)$ have
uniform geometry,  so the norm of the Sobolev embedding
$W^{1,p}\hookrightarrow L^\infty$  on $C_z(\delta)$ is bounded
independent of $z$.  Because $\rho$ is bounded, so we have
$\|\xi\|_\infty \leq c(p)\|\xi\|_{1,p,\la}$ on $C_z(\delta)$.
Also, on each component of $B_z(\delta)$, the proof of
Lemma~\ref{diamlemma} shows that $\mbox{osc} |\xi| \leq
c\|\nabla\xi\|_{p,\la}$.  The lemma follows by writing $\xi=\zeta+\bar{\xi}$ as in (\ref{sumnorm}) and noting that, because $\bar{E}_f$ is finite dimensional, there is a constant so that $\|\nabla\bar{\xi}\|_{1,p,\la} \leq c\|\bar{\xi}\|$ for all $\bar{\xi}\in \bar{E}_f$. \qed

\bigskip

Lemma~\ref{supbound} implies a bound on the norm of the
linearization $L_f$.   The next two lemmas give more specific -- and crucial -- information.  They show 
that, for $p>2$, $L_f$ is uniformly bounded and continuous in a  neighborhood of the
space of $J_\a$-holomorphic maps in the  $\la_p$-topology of Definition~\ref{defMapE}.

\medskip

\begin{lemma}[Uniform boundedness of $L_f$]
\label{Lfbounds} For each $p>2$ there is a constant $c=c(p)$ and a
neighborhood  $\N$ of the space of $J_\a$-holomorphic maps in the
$\la_p$-topology such that for each  $f:C\to N_D$ in $\N$ and each
$\xi\in \E_f$
$$
\|L_f\xi\|_{p,\la} \leq c(1+E) \, \|\xi\|_{1,p,\la}
$$
where $E$ is the constant of Lemma~\ref{compactnesstheorem}.
\end{lemma}

\pf   From (\ref{Lfformula}) we have the pointwise bound $\frac12
|L_f\xi |\leq (1+|\alpha|) |\nabla\xi| +|\nabla \a| |\xi||df|$.  By
Lemma~\ref{ILLemma} the image of each $J_\a$-holomorphic map lies in
$D$. It then follows from  Lemma~\ref{Hdistlemma} that there is a
neighborhood $\N$ of the space of $J_\a$-holomorphic maps in
$\la_p$-topology such that the images of all maps in $\N$ lie in a
compact set in $N_D$ where $|\alpha|+|\nabla \alpha|<c_1$ for some
$c_1$.  Integrating using norms (\ref{1pnorm}) and (\ref{pnorm}), we
have
$$
\|L_f\xi\|_{p,\la} \leq (2+2c_1) \|\xi\|_{1,p,\la} + 2 c_1
\|\xi\|_\infty\cdot \|df\|_{p,\la}
$$
with $\|\xi\|_\infty$  bounded as in (\ref{supboundeq}) and
$\|df\|_{p,\la} \leq E$ by Lemma~\ref{compactnesstheorem}. The lemma
follows. \qed

\bigskip

\subsection{Estimates at a fixed map}

We next establish elliptic estimates for $L_f$ when $f$ is a fixed $J_\alpha$-holomorphic map.  Such estimates are standard for unweighted Sobolev spaces, but we need them for   the weighted spaces $\E_f$ and $\F_f$.  We start by considering $L_f$ acting on functions on $\E_f^0$ on a neck region $B(\de)$.

\begin{lemma}[Interior estimate]
\label{interiorlemma}
For $p>2$  and each $J_\alpha$-holomorphic map $f$,  there are constants $\delta_1$ and $c$, depending on $p$ and $f$,  such that whenever  $\delta<\de_1$ and $L_f\xi\in\F_f$ with $\|\xi\|_{p,\la}<\infty$,  then $\xi\in\E_f^0$ on  each neck region $B(\de)$,  and 
\begin{equation}
\label{interiorestimate}
\|\xi\|_{1,p,\la}\leq c\, \big(\|L_f\xi\|_{p,\la} + \|\xi\|_{p,\la}\big)
 \qquad \forall \xi\in \E_f^0 \mbox{  on $B(\delta)$}.
\end{equation}
\end{lemma}

\pf  By the Image Localization Lemma the image of $f$  lies in the set where $\alpha=0$,  and by Lemma~\ref{diamlemma} we can
choose $\delta$ small enough that the image of  the neighborhood $B(\delta)$ of each node lies in a holomorphic coordinate chart for the bundle $N\to N_D$.   In such a  chart, we can regard $f$ and  any $\zeta\in \E^0_f $ as complex-valued functions on $B(\delta)$ or, equivalently, on $T(\delta)$ via (\ref{defphi}).   Using formula (\ref{LRformula}) we can then write   
\begin{equation}
\label{Ldel}
L_f\zeta = \del\zeta + \Gamma\zeta + R_\alpha\zeta
\end{equation}
where $\del$ is the standard $\del$-operator on functions, $\Gamma$ is a term built from the Christoffel symbols satisfying
 $|\Gamma|\leq c_2\,\diam(f(B(\delta))) |df|\leq c_3|df|$, and where $|R_\alpha|\leq c_4\,|df|$.

 For each integer $n$ let $A(n)\subset T$ be the annulus with $n\leq t \leq n+1$ and let $A^+(n)$ be the annulus with $n-1\leq t \leq n+2$.  As in Appendix C of \cite{MS},  standard elliptic theory shows that if $\zeta\in L^p$   and $L_f\zeta\in L^p$ on $A_n^+$, then $\zeta\in W^{1,p}$ on $A_n$ and
  \begin{equation}
\label{firstelliptic}
 \|\zeta\|_{1, p;A(n)} \leq c_5\,\left(\|\del\zeta\|_{p;A^+(n)} +\|\zeta\|_{p;A^+(n)}\right)
\end{equation}
 where these are unweighted $L^p$  norms.  The constant  $c_5$ is independent of $n$ because the $A^+(n)$ are isometric to one another. 
 Summing on $n$, noting that each point on $T(\delta)$ lies in three of the $A^+(n)$, and inserting (\ref{Ldel}) with the associated bounds on $\Gamma$ and $R_\alpha$ yields
  \begin{equation}
\label{interioresteq3}
 \| \zeta\|_{1, p; T(\delta)} \leq 3c_6\, \big(\|L_f\zeta \|_{p; T(\delta)} +\| \zeta\|_{p; T(\delta)}\ + c_7\,\| | \zeta| df\|_{p, T(\delta)} \big).
\end{equation}

Now take $\zeta=\rho^{-7\la/p}\xi$ and define operators $\nabla^\la=\rho^{7\la/p}\nabla\rho^{-7\la/p}$ and $L_f^\la=\rho^{7\la/p}L_f\rho^{-7\la/p}$.  Then (\ref{interioresteq3}) becomes
  \begin{equation}
\label{interioresteq4}
 \int_{T(\delta)} \rho^{-7\la}  \bigg(|\nabla^\la\xi|^p + |\xi|^p\bigg) \  \leq \ c_8\,  \int_{T(\delta)} \rho^{-7\la}  \bigg( |L^\la_f\xi|^p + |\xi|^p + |df|^p \|\xi\|^p_\infty \bigg).
\end{equation}
But $\nabla^\la\xi=\nabla\xi-\frac{7\la}{p}d\log\rho\cdot \xi$ and $L_f^\la\xi=L_f\xi- \frac{7\la}{p}(d\log\rho\cdot \xi)^{0,1}$ and the definition of $\rho^2$ given after (\ref{defphi}) implies that $|d\log\rho|=|\tanh(2t)|\leq 1$.  Hence we can replace $\nabla^\la$ by $\nabla$ and $L_f^\la$ by $L_f$ in (\ref{interioresteq4}), absorbing the difference into the $|\xi|^p$ term.

After rewriting (\ref{1pnormB}) as  integrals over $B(\delta)$ via (\ref{alphanorm2}), we can   bound  $\|df\|_{p,\la}$ by  (\ref{aboundonM}) and $\|\xi\|_\infty$ by  (\ref{supboundeq}) to obtain  
\begin{equation}
\label{interioresteq5}
 \|\xi\|_{1,p,\la; B(\delta)} \  \leq \   c_{11}\left(\|L_f\xi\|_{p, \la; B(\delta)} +\|\xi\|_{p, \la; B(\delta)}\right) \ +c_{12}\delta^{\frac{1}{4p}} \|\xi\|_{1,p,\la} .
\end{equation}
These constants are independent of $\delta$, so we can choose $\delta$ small enough that the last term can be absorbed into the left-hand side, giving the stated  inequality (\ref{interiorestimate}).

Finally, if  the weighted $L^p$ norms of $\xi=\rho^{7\la/p}\zeta$ and $L_f\xi$ are finite,  then $\zeta$ and $L_f\zeta$ are in $L^p$. Hence, as noted above (\ref{firstelliptic}), $\zeta$ is in $W^{1,p}$ on $T(\de)$.  The conversion back into integrals of $\xi$ done above then shows that the weighted $(1,p)$ norm of $xi$ is finite, so $\xi\in\E^0_f$.
\qed

\bigskip

\begin{prop}[Elliptic estimate]
\label{basicellipticlemma}
For $p>2$ each $J_\alpha$-holomorphic map $(C,f)$ there is a constant $c$, depending on $p$ and $f$,  such that whenever $\xi\in L^2$ and $L_f\xi\in\F_f$ then $\xi\in\E_f$ and 
\begin{equation}
\label{basicelliptic}
\|\xi\|_{1,p,\la}\leq c\, \big(\|L_f\xi\|_{p,\la} + \|\xi\|_{p,\la}\big)
 \qquad \forall \xi\in \E_f.
\end{equation}
\end{prop}

\pf  For any $\delta>0$, the weighting function $\rho$ is bounded above and below on $C(\delta)$, so the weighted norms (\ref{1pnorm}) and (\ref{pnorm})  are equivalent to the usual Sobolev norms. Again, elliptic theory  shows that if $\xi\in L^2$   and $L_f\xi\in L^p$, then $\xi\in W^{1,p}$ on $C(\delta)$  and
\begin{equation}
\label{basicelliptic2}
\|\xi\|_{1,p, \la; C(\delta)}\leq c_1(p,\delta)\, \big(\|L_f\xi\|_{p, \la; C(\delta)} + \|\xi\|_{p, \la; C(\delta)}\big).
\end{equation}
 In particular, this proves the lemma  whenever the domain $C$ of $f$ is smooth.

When $C$ is  nodal, let  $\delta_1$ be  the constant of  Lemma~\ref{interiorlemma} and let $\beta_{\de_1}$  be the corresponding cutoff function defined before  (\ref{defxibar}).  Given $\xi$,  write 
$$
\xi=\xi_1+\zeta+\bar{\xi}
$$
 where   $\xi_1=(1-\beta_{\de_1})\xi$ and $\zeta+\bar{\xi}$ is the decomposition of $\beta_{\de_1}\xi$ as in (\ref{defxibar}) for this $\de_1$.  Then
\begin{enumerate}
\item[(i)]  $\xi_1$ satisfies (\ref{basicelliptic2}) and $\|L_f\xi_1\|_{p,\la}=\|(1-\beta)L_f\xi +d\beta\cdot\xi\|_p \leq \|L_f\xi\|_{p,\la}+c_2(\delta)\|\xi\|_{p,\la}$.
\item[(ii)] Similarly, $\zeta$ satisfies the estimate of Lemma~\ref{interiorlemma},  and $\|L_f\zeta\|_{p,\la} \leq \|L_f\xi\|_{p,\la}+c_2(\delta)\|\xi\|_{p,\la}$.
\item[(iii)]   $\bar{E}_f$ is finite-dimensional and $L_f$ is linear and injective by the Vanishing Theorem~\ref{Vtheorem}.  Hence there is a constant such that $\|\bar{\xi}\|_{1,p, \la}\leq c_3\|L_f\bar{\xi}\|_{p,\la}$ for all $\bar{\xi}\in \bar{E}_f$.  Furthermore, 
$$
 \|L_f\bar{\xi}\|_{p,\la} \,=\,  \|L_f\xi-L_f(\xi_1+\zeta)\|_{p,\la}
  \,\leq\,  c_4\,\big(\|L_f\xi\|_{p,\la} + \|L_f\xi_1\|_{p,\la} +\|\zeta\|_{p,\la}\big).
  $$
\end{enumerate}
Altogether, we have
\begin{eqnarray*}
\|\xi\|_{1,p, \la} & \leq & c_5\left(\|\xi_1\|_{1,p, \la} + \|\zeta\|_{1,p, \la} +\|\bar{\xi}\|_{1,p, la}\right) \\
& \leq & c_6\big( \|L_f\xi_1\|_{1,p, \la} + \|L_f\zeta\|_{1,p, \la} +\|L_f\bar{\xi}\|_{1,p, \la} +\|\xi\|_{1,p, \la} \big)\\
& \leq & c_7\big( \|L_f\xi\|_{1,p, \la}  +\|\xi\|_{1,p, \la} \big). \hspace{5cm} \mbox{\qed}
\end{eqnarray*}

\bigskip

\begin{lemma}
\label{firstestimate} For each $p> 2$ and each
$J_\alpha$-holomorphic map $f:C\to N_D$ there is a constant $c(p,
f)$ such that  
\begin{equation}
\label{firstestimateEq} \|\xi\|_{1,p, \la}\ \leq \ c(p,f)\, \|L_f\xi\|_{p, \la}
\qquad\forall  \xi\in \E_f.
\end{equation}
\end{lemma}

\pf  Otherwise there is a sequence $\xi_k\in \E_f$  such that
$\|\xi_k\|_{1,p, \la} \geq k\|L_f\xi_k\|_{p, \la}$; dividing by
$\|\xi_k\|_{1,p, \la}$ we can assume that $\|\xi_k\|_{1,p, \la}=1$ for all
$k$.   Then $\|L_f\xi_k\|_{p, \la} \to 0$ and, by the weak compactness of the unit ball in a Banach space, we can pass to a subsequence that converges weakly to some $\xi_0\in \E_f$.  Furthermore, using Holder's inequality, (\ref{supboundeq}), and the fact that the weighting function $\rho$ of (\ref{defrho}) is bounded, we have
\begin{eqnarray}
\label{5.L12}
\int_C |\nabla\xi|^2 + |\xi|^2 & \leq & \big(\mbox{vol(C)}\big)^{\frac{2p-2}{p}}\left(\int_C |\nabla\xi|^p\right)^{\frac{2}{p}} + c_1 \|\xi\|^2_\infty \nonumber\\
& \leq & c_2 \left(\int_C \rho^{-\la}|\nabla\xi|^p\right)^{\frac{2}{p}}  + c_3\|\xi\|^2_{1,p, \la}\nonumber\\[.1cm]
& \leq &  c_4\|\xi\|^2_{1,p, \la}.
\end{eqnarray}
Thus $\{\xi_k\}$ is bounded in the unweighted $W^{1,2}$ norm.  By the compactness of the
embedding $W^{1,2}\hookrightarrow L^2$ on the fixed curve $C$,
there is a $L^2$-convergent subsequence $\xi_k\to \xi_0$.  Applying
the elliptic estimate (\ref{basicelliptic}) to $\xi_k-\xi_\ell$ then shows that  $\xi_k\to
\xi_0$ in the $(1,p, \la)$ norm.  But then $\|\xi_0\|_{1,p, \la}=1$, so $\xi_0$ is a
non-zero solution of  $L_f\xi_0=0$,  contradicting the Vanishing Theorem~\ref{Vtheorem}. \qed

\bigskip

\subsection{Estimates on $L_f$ uniform on the space of maps}

The constant in inequality  (\ref{firstestimateEq}) depends on
the   $J_\alpha$-holomorphic map $f$.  Our last set of analysis results show that  the same inequality holds with a constant that is uniform for all $f$ in a neighborhood of the space of  $J_\alpha$-holomorphic maps.   Again the issue is uniformity around  the maps with nodal domains.

\medskip

\begin{lemma}[Continuity of $L_f$]
\label{continuitylemma}
Fix a $J_\alpha$-holomorphic map $(C_0,f)$.  For each $p> 2$ and each $\ep>0$ there is a  neighborhood $\N_f(p,\ep)$  of $(C_0,f)$ in the space of maps such that, for any $g\in \N_f(p,\ep)$  there are Banach space maps $P:\E_g\to \E_f$  and $Q:\F_g\to \F_f$ with
\begin{equation}
\label{Lf-Lg}
\|L_g\xi-(Q^{-1}L_fP)\xi\|_{p,\la} \ \leq\ \ep \|\xi\|_{1,p, \la} \qquad \forall \xi\in \E_g.
\end{equation}
Moreover,  $A=P, Q$ both  satisfying  $(1-\ep)\|\xi\|_{1,p, \la}\leq \|A\xi\|_{1,p, \la}\leq (1+\ep)\|\xi\|_{1,p, \la}$.
\end{lemma}

\pf   We will compare $L_f$ with  $L_g$ using the explicit  formula (\ref{Lfformula}).  First note that the proof of Lemma~\ref{Lfbounds} produced a
neighborhood $\N$ of $f$ so that all maps  in $\N$ have images in a compact subset of $N_D$.  For points $x,y$ in that compact subset the smooth tensor fields $J$ and $K_\alpha$ satisfy
\begin{equation}
\label{JKbound}
|J_x-J_y| + |(K_\alpha)_x - (K_\alpha)_y| + |(\nabla K_\alpha)_x - (\nabla K_\alpha)_y| \ \leq\ c_1\,\dist(x,y)
\end{equation}
 and $|J|+|K_\alpha|+|\nabla K_\alpha|<c_2$ with constants independent of $x$ and $y$.   Furthermore, Definition~\ref{defMapE} implies that there is  a neighborhood $\N_f\subset \N$ so that the domains of all maps in $\N_f$ lie in a single chart  in the universal curve of the form (\ref{firsttrivialization}) or (\ref{secondtrivialization}).  The domain of each $(C, g)\in\N_f$ can then be written as  the union of $C(\delta)$ and $B_i(\delta)$ as in (\ref{defCd}), and the chart gives an identification  $C(\delta)= C_0(\delta)$ that is $C^1$ close to an isometry.
 
 When $C_0$ is smooth the chart (\ref{firsttrivialization})  identifies the domains $C=C_0$, so we can regard all maps in $\N_f$ as maps $g:C_0\to N_D$.  Lemma~\ref{diamlemma} and \ref{Hdistlemma} then imply that 
 \begin{equation}
\label{JKbound2}
 \sup_{x\in C_0} \dist(f(x), g(x))\leq c_3 \vv f-g\vv_{1,p, \la}
 \qquad\forall g\in\N_f
\end{equation}
where   $ \vv f-g\vv_{1,p}$ denotes the sum of the distance between $C$ and $C_0$ in $\bM_{g,n}$ and $\|f-g\|_{1,p, \la}$ defined by (\ref{alphanorm}).

 Hence, after possibly making $\N_f$ smaller, there is a unique geodesic from $f(x)$ to $g(x)$ for each $x$, and parallel transport along this geodesic induces linear maps $P:\E_g\to \E_f$  and $Q:\F_g\to \F_f$.   this gives two connections on $g^*N$ over $C_0$: the pullback connection $g^*\nabla$ of the connection $\nabla$ on $N$, and the connection $Q^{-1}(f^*\nabla)P$ obtained from the pullback $f^*\nabla$.   The difference   $\Gamma=Q^{-1}\nabla^fP - \nabla^g$ is an $\mbox{End}(N)$-valued 1-form with $|\Gamma(x)|\leq c_4\, \dist(f(x), g(x))$.   Referring to the formula (\ref{Lfformula})  and using (\ref{JKbound}) and (\ref{JKbound2}), we then have the pointwise bound
\begin{eqnarray*}
 \left|L_g\xi-(Q^{-1}L_fP)\xi\right|  & \leq &  c\Big(       
 |J_f-J_g| |\nabla\xi| + |\Gamma(\xi)| +|(\nabla K_\alpha)_f-(\nabla K_\alpha)_g|\, |\xi|\, |df| \\
& \ & \hspace{4.5cm}  + |\nabla K_\alpha|\, |\xi|\, |df-dg| + |\alpha| \, |\nabla\xi|  \Big)\\[.2cm]
 & \leq &  c\Big( \vv f-g\vv_{1,p}\left(|\nabla\xi| +  |\xi|\, |df|\right) +  |\xi|\, |df-dg| + |\alpha| \, |\nabla\xi|\Big).
 \end{eqnarray*}
 We also have the sup bound (\ref{supboundeq}) and, because $\alpha$ vanishes along the image of $f$, $|\alpha|\leq c_5 \|f-g\|_\infty$.  Hence the above inequality gives
  \begin{equation}
\label{JKbound3}
  \| L_g\xi-(Q^{-1}L_fP)\xi\|_{p, \la}  \  \leq \   c_6 \|\xi\|_{1,p, \la}\cdot \vv f-g \vv_{1,p}\, (1+\|df\|_{p, \la}).
\end{equation}
  By shrinking the neighborhood $\N_f$, we can make the last two factors as small as desired, giving  (\ref{Lf-Lg}).
  
  When $C_0$ is nodal, the above analysis applies on $C(\delta)$ for each $\delta>0$, so it suffices to prove a similar bound for the neighborhood $B(\delta)$ of a node for some small $\delta$.  In fact,  by Lemma~\ref{diamlemma}, we can choose $\delta$ so that  the image $f(B(\delta))$ lies in a region  that can be covered by  a geodesically convex local holomorphic chart for the holomorphic bundle $\pi^*N\to N_D$ (cf. (\ref{split}).   We can thus regard $f$ and $\xi$ as complex-valued functions.  Again, parallel transport defines isometries $P, Q$ that are $C^1$ close to the identity in coordinates, and under which the 1-form  $\Gamma=f^*\nabla-g^*\nabla$ is pointwise bounded by   $|\Gamma|(x)\leq c_7\, \dist(f(x), g(x))$.
 
      Now  identifying $B(\delta)$ with the cylinder $T$ as in (\ref {defphi}), we have $L_g\xi=\nabla_t\xi\, dt + J\nabla_\theta\xi\,d\theta$, and similarly for $Q^{-1}L_fP\xi$.   We can then bound $ \left|L_g\xi-(Q^{-1}L_fP)\xi\right|$ pointwise exactly as above, and integrate to obtain (\ref{JKbound3}) on $B(\delta)$.
      
      Finally, in coordinates on both $C(\delta)$ and $B(\delta)$ the parallel transport operators $A=P, Q$ both satisfy  the pointwise bound $|A\xi-\xi|_{C^1}\leq c_8 \|f-g\|_\infty \leq c_9\vv f-g \vv_{1,p}$.  Integrating (and possibly shrinking $\N_f(p, \ep)$) then gives the inequalities on the last line of the lemma for all $g\in\N_f(p, \ep)$.
 \qed

\begin{theorem}
\label{MainKernelZeroTheorem}
  For each $p> 2$ there is a constant $c(p)$ and a neighborhood $\N$ of the space of stable $J_\alpha$-holomorphic maps in $\Map_{g,n}(N_D,d)$ 
such that, for every $f\in \N$, 
\begin{equation}
\label{finalindex}
L_f:\E_f\to \F_f
\end{equation}
is a uniformly bounded  Fredholm map with $\ind L_f=-2\beta$ as in  (\ref{defbeta})  and with  
\begin{equation}
\label{MainKernelZeroTheoremeq}
\|\xi \|_{1,p, \la}\ \leq c(p)\, \|L_f\xi \|_{p, \la}
\qquad \mbox{$\forall \xi\in \E_f$}.
\end{equation}
\end{theorem}

\pf  For each stable $J_\alpha$-holomorphic map $f$, let $c(p,f)$ be the associated constant of Lemma~\ref{firstestimate} and let 
$\N_f$ be  the neighborhood $\N_f(p, \ep)$ provided by Lemma~\ref{continuitylemma} with $\ep=1/(2+2c(p,f))$.   Applying Lemmas~\ref{continuitylemma} and (\ref{firstestimate}), we see that  each $\xi\in E_g$ satisfies
\begin{eqnarray*}
\label{5.12}
\|\xi\|_{1,p,\la}\, =\, \|P\xi \|_{1,p, \la}& \leq & c(p,f)\|L_fP\xi \|_{p, \la} \\
& \leq & c(p,f)\|Q^{-1} L_g\xi\|_{p, \la}\, +\, c(p,f)\|Q^{-1}L_fP\xi-L_g\xi\|_{p,\la}.
\end{eqnarray*}
But by Lemma~\ref{continuitylemma} $\|Q^{-1}\|\leq(1-\ep)^{-1}< 2$ and the last term is bounded by $\frac12 \|\xi\|_{1,p, \la}$.   Hence
$$
\|\xi\|_{1,p, \la}\ \leq \  3c(p,f)\|L_g\xi\|_{p, \la}
\qquad \mbox{$\forall \xi\in \E_f$}
$$
whenever $g\in \N_f$.  The sets $\{\N_f\}$ obtained in this  manner cover the space of stable maps, so (\ref{MainKernelZeroTheoremeq}) follows by the Compactness Theorem~\ref{compactnesstheorem}. 

We know from Theorem~\ref{Lfbounds} that $L_f$ is uniformly  bounded for $f\in\N$, and  $\ker L_f=0$ by  (\ref{MainKernelZeroTheoremeq}).  Inequality (\ref{MainKernelZeroTheoremeq})  also implies that the range of $L_f$ is closed: if $L_f\xi_k\to \eta$ then (\ref{MainKernelZeroTheoremeq}) (applied to $\xi_k-\xi_\ell$) shows that $\{\xi_k\}$ is Cauchy, so $\xi_k\to\xi_0$ in $W^{1,p}$ with $L_f\xi_0=\lim L_f\xi_k=\eta$, so $\eta\in\im\,L_f$.   The proof is completed by noting that $\dim \cok L_f=2\beta$ by Lemma~\ref{5.lastlemma} below.
\qed

\begin{lemma}
\label{5.lastlemma}
For each $f$ in the neighborhood $\N$ of Theorem~\ref{MainKernelZeroTheorem}, we can choose a $2\beta$-dimensional subspace $V_f\subset \F_f$, consisting of smooth forms that vanish in a neighborhood of the nodes, that is transverse to the image of $L_f:\E_f\to \F_f
$.  Hence $\dim \cok L_f=2\beta$.
\end{lemma}

\pf  First let $E^{1,p}_f$ and $F^p_f$ be, respectively  the completions of $E_f$ and $F_f$ defined by (\ref{DefEF}) in the usual, unweighted $W^{1,p}$  and $W^{0,p}$ norms.    By  Lemma~12.2 of \cite{FO}  (see also Theorem~C.1.10 of \cite{MS}) the linearization (\ref{Lbounded}) extends to a bounded Fredholm map
$$
L_f: E^{1,p}_f\to F^p_f
$$
whose index is $-2\beta$ and whose kernel vanishes by Theorem~\ref{Vtheorem}.  Hence we can choose  linearly
independent (0,1) forms $\{v_j\,|\,j=1,\dots,2\b\}$ in $F^p_f$ so that
$V_f=\mbox{span}\{v_j\}$ is transverse to $\im\, L_f$.  We can assume that each $v_i$ is smooth and vanishes in a  neighborhood of all nodes because such forms are dense  in $L^p$.  Then each $\eta\in F^p_f$ can be uniquely written as
\begin{equation}
\label{eta=Lv}
 \eta=L_f\xi + v 
 \end{equation}
for some $\xi\in E^{1,p}_f$ and $v\in V_f$.  We will show that the same decomposition holds for the weighted Sobolev spaces $\E_f$ and $\F_f$  defined by the weighted norms (\ref{1pnorm}) and (\ref{pnorm})

 Because  $\rho$ is bounded,  each $\eta\in \F_f$ satisfies
$$
\int_C |\eta|^p\ \leq\  c \int_C \rho^{-\la} |\eta|^p.
$$
Thus $\eta$ is in the unweighted space $F^p_f$, so we can write $\eta=L_f\xi+v$ as above.  Furthermore, 
 $\rho^{-1}$ is bounded on the support of each $v_i$, so $V_f\subset\F_f$.   We can then apply Proposition~\ref{basicellipticlemma} to
 $L_F\xi=\eta-v \in \F_f$ to conclude that $\xi\in \E_f$.  Thus  each $\eta \in \F_f$ can be written as the sum of an element  in $L_f\E_f$  and a $v\in V_f$.  Moreover, this decomposition is unique:  if  $L_f\xi=v$ for some $\xi\in\E_f$ and $v\in V_f$ then $\xi\in E_f^{1,p}$ by (\ref{5.L12}), contradicting the uniqueness of (\ref{eta=Lv}).  We conclude that $\dim \cok L_f=\dim V_f =2\beta$.
\qed

 \bigskip
 
Theorem~\ref{MainKernelZeroTheorem} and the other results of this section apply to the linearization $L_f$ at a single map $f$.   We will next examine how $L_f$ as $f$ varies in the space of maps.

\vspace{1cm}

\setcounter{equation}{0}
\section{Obstruction bundle $\Ob$}
\label{section6}
\medskip

Abstractly, $\Ob$ is the topological vector bundle whose fiber at a $J_\a$-holomorphic map $f:C\to
N_D$ is $\cok L_f$.  Thus it is defined by the exact sequence
$$
0\to \E \overset{L}\longrightarrow \F\overset{\rho} \longrightarrow
\Ob \to 0
$$
of topological vector bundles over $\Map(X)$. The goal of this
section is to give a concrete realization of $\Ob$ and show that  it is a locally trivial bundle on the space of
maps. For this, we will split the sequence
\begin{equation}
\label{shortsequence} 0\to \E_f\overset{L_f}\longrightarrow
\F_f\overset{\rho} \longrightarrow \cok L_f \to 0
\end{equation}
 by regarding $\cok L_f$ as a subspace of $\F_f$ for each map $f$.  The key issue is constructing a splitting that is locally trivial around maps with nodal domains.

\medskip
Consider the space  $\Map_{g,n}(N_D,d)$  of  maps  into  the open
complex surface $N_D$ that represent the class $d[D]$.  The following theorem is the main application of the analysis done
in Sections~2 and 5.

\begin{theorem}
\label{constuctionofOb} There is a  locally trivial real vector
bundle $\Ob$  defined in a $\la_p$-topology neighborhood  $\N$ of
$\CM_{g,n}(D,d)$    in $\Map_{g,n}(N_D,d)$ such
that
the fiber of $\Ob$ at each $f\in {\cal N}$  lies in  $F_f$, is
transversal to the image of  $L_f$, and has the same dimension as
$\cok L_f$.   This bundle $\Ob$ has a canonical orientation.
\end{theorem}

The remainder of this section proves Theorem~\ref{constuctionofOb}
in three steps. The first step adopts the method of Fukaya and Ono
\cite{FO}.
Step~2, which shows that the obstruction bundle has constant rank,
is the key feature that is special to our situation. Step~3 shows
how the local trivializations define a global bundle. \vskip.5cm

\noindent{\bf STEP 1 --- The fiber  $\Ob_f$ at  a map $f:C\to N_D$.}\ \
For a fixed $f\in \CM_{g,n}(D,d)$ the subspace $V_f$ of Lemma~\ref{5.lastlemma} specifies a (non-canonical) splitting of 
(\ref{shortsequence}.   Again, write $V_f$ as the space of  linearly
independent (0,1) forms $\{v_j\,|\,j=1,\dots,2\b\}$ in $F_f$.  Also  consider the subspace
\begin{equation}\label{01-form}
{\cal P}_0= \Omega_{0}^{0,1}(\bU\boxtimes N) 
\end{equation}
of the space (\ref{defP}) of Ruan-Tian perturbations that have values in the subspace $N\subset TN_D$ and
vanish near all nodes. For each $v$  in the space (\ref{01-form}) the restriction of $v$ to the graph  of $f$, defined as in (\ref{restriction}), is an element $v_f\in \Omega^{0, 1}_{0}(C,f^*N)$ that vanishes  in all of the regions $B(2\delta_0)$ around the  nodes.  The following lemma shows that, conversely,  all such elements of  $\Omega^{0, 1}_{0}(C,f^*N)$ arise as restrictions of elements of ${\cal P}_0$.

\begin{lemma}
\label{extensionLemma}
Given $v\in  \Omega^{0, 1}_{0}(C,f^*N)$, there exists
$\hat{v}$ in the space (\ref{01-form}) with $\hat{v}_f=v$.
\end{lemma}

\pf    Let $C(\delta_0)$ be the part of $C$ away from the nodes where $\rho\geq
\delta_0$. The graph $G_f=\{(x,f(x))\,|\,x\in C(\d_0)\}$ is a
submanifold of $C(\d_0)\times N_D$ that has a $2\ep$-tubular
neighborhood $W$  consisting of points $(x,y)$ such that
there is a unique minimal geodesic in $N_D$ from $y$ to $f(x)$.
Parallel transporting the fibers of $N$ along these geodesics
trivializes $N$ over $W$.  Fix a cutoff function $\beta_W\in
C^\infty(C(\d_0)\times N_D)$ with support on $W$ and with
$\beta_W\equiv 1$ on the $\ep$-tubular neighborhood of $G_f$.  Then,
regarding $v$ as a form on $G_f$, we can extend to $W$ by parallel
transport, multiply by $\beta_W$, and extend by zero to obtain
a section $\beta_W v$ of  $T^*C(\d_0)\boxtimes N$ over
$C(\d_0)\times N_D$.

Now  fix  a local trivialization $\phi:C(\delta_0)\times V\to U$ of the universal curve $\bU$  around $C$ as in (\ref{secondtrivialization}) and write $\phi^{-1}=
(\tau_1, \tau_2)$.  Choose
a smooth cutoff function $\beta_U$ on $\bU$ with support on
$\pi^{-1}(U)$ and with $\beta_U\equiv 1$ on a smaller neighborhood
$\pi^{-1}(U')$ of $C$.    The section $\tilde{v}=\beta_U
\tau_1^*(\beta_Wv)$ then extends by zero to  a section of
$T^*\bU\boxtimes N$ over   $\bU\times N_D$. The $Jj$-antilinear
component $\hat{v}=\tfrac12(\tilde{v}+J\circ \tilde{v}\circ j)$ is
the desired (0, 1) form. \qed

\medskip

\bigskip

\noindent{\bf STEP 2 --- The local trivialization of the cokernel
bundle.}\ \  Fix  $f\in\CM_{g,n}(D,d)$ and $V_f=\mbox{span}\{v_j\}$   as in Step~1.  In light of
Lemma~\ref{extensionLemma}, we can  extend the forms in $V_f$ to obtain a vector space
$$
\hat{V}_{f}\,=\,\mbox{span}\{\,\hat{v}_j\,\}\ \subset\
{\cal P}_0.
$$
Given   a neighborhood $\N_f$ of $f$ in the space of maps, we can regard $\hat{V}_{f}$ as a trivial vector bundle  $\hat{V}_{f} \times {\cal
N}$ over ${\cal N}$.  For each map $g\in\N_f$ the composition of the
restriction (\ref{restriction}) and the projection $\rho$ in (\ref{shortsequence})  
gives a linear map $R$ defined by  $R(\hat{v},g)=\rho(\hat{v}_g)$  and a diagram
\begin{equation}
\label{projtocoker}
\begin{xy}
\xymatrix@!0 @R=1.2cm @C=1.3cm{
\hat{V}_{f} \times \N_f  \ar[dr]    \ar[rr]^-R  & & \cok L  \ar[dl]  \\
  &  \N_f   &  }
 \end{xy}
\end{equation}
where $\cok L$ is the vector bundle whose fiber at $g$ is $\cok L_g$.

\medskip

\begin{lemma}
\label{7.analysis} The  vector bundle $\cok L$ is locally trivial:  each $f\in\CM_{g,n}(D,d)$ has a neighborhood $\N_f$ in the space of maps for  which (\ref{projtocoker}) is a vector bundle isomorphism.
\end{lemma}

\pf  By Theorem~\ref{MainKernelZeroTheorem} there is a neighborhood
${\cal N}$ of $\CM_{g,n}(D,d)$    such that for every
$g\in {\cal N}$ we have  $\ker L_g=0$ and $\dim \cok L_g= 2\beta$.  In particular, for each $g\in \N$, $\cok L_g$ has the same dimension as $\hat{V}_{f}$.   It suffices to find a neighborhood $\N_f$ on which $R$ is injective on fibers or, equivalently, on which $\hat{V}_{f} \cap \mbox{im\,$L_g$}=0$.

Let $S(V)$ denote the unit sphere in the (finite-dimensional) space $\hat{V}_{f} $ defined by the condition $\|v\|_\infty+\|\nabla v\|_\infty=1$.  Because $S(V)$ is compact  and $\mbox{im}\,L_g$ is closed, the infimum
\begin{equation}
\label{defepg}
\ep_g\,=\,\inf \big\{ \|L_g\xi-\hat{v}_g\|_{p, \la}\, :\, \xi\in \E_g \mbox{ and } v\in S(V)\big\}
\end{equation}
is realized for each $g$ and is equal to 0 if and only if $\hat{V}_{f} \cap \mbox{im\,$L_g$}\neq 0$.  Thus we have $\ep_f >0$ (by our choice $\hat{V}_{f}$), and must show that $\ep_g>0$ for all nearby $g$.

    Lemma~\ref{continuitylemma} shows that there is a neighborhood $\N_f(\ep)$  such that for each $g\in \N_f(\ep)$ satisfies $\|f-g\|_{1,p, \la}<\ep$ and we have the inequalities   (\ref{Lf-Lg}) and $\|Q^{-1}\eta\|_{p, \la}\geq (1-\ep)\|\eta\|_{p, \la}$ for all $\eta\in\E_f$. Hence for each  $(\xi, v)\in \E_g\times S(V)$ we have
\begin{eqnarray*}
 \|L_g\xi -\hat{v}_g\|_{p, \la} & \geq & \|Q^{-1}(L_fP\xi -\hat{v}_f)\|_{p, \la}\ - \|L_g\xi-Q^{-1}L_fP\xi\|_{p, \la} - \|Q^{-1}\hat{v}_f -\hat{v}_g\|_{p, \la} \\[.2cm]
& \geq & (1-\ep) \|L_f\xi- Q\hat{v}_g\|_{p, \la}  - \ep\, \|\xi \|_{1,p} - \|Q^{-1}\hat{v}_f -\hat{v}_g\|_{p, \la}.
\end{eqnarray*}
with $ \|L_f\xi- \hat{v}_f\|_p\geq \ep_f>0$ by the  Definition (\ref{defepg}).  Furthermore, 
\begin{itemize}
\item using Theorem~\ref{MainKernelZeroTheorem} we have $\|\xi\|_{1,p, \la} \leq c_1\|L_g\xi\|_{p, \la} = \| \hat{v}_g\|_{p, \la}  \leq c_2\|\hat{v}_g\|_\infty \leq c_3$, 
\item by the definitions  of  $Q$ and $S(V)$, $\|Q^{-1}\hat{v}_f -\hat{v}_g\|_{p, \la}$ is dominated by $ 2\|\nabla v||_\infty \|f-g\|_{\infty}  \leq c_4 \|f-g\|_{1,p, \la}$.
\end{itemize}
Taking $\ep$ sufficiently small, we see that $\ep_g>0$ for all $g\in \N_f(\ep)$, as required.
\qed

\bigskip

\noindent{\bf  STEP 3 ---  Definition of $\Ob$.}\ \ The procedure
of Step~2 can be applied to each $J_\a$-holomorphic map $f$ to
obtain a pair $(\hat{V}_f,\N_f)$. By compactness, we can choose a
finite set $\{f_i\}$ of $J_\a$-holomorphic maps so that
$\{\N_{f_i}\}$ covers $\CM_{g,n}(D,d)$. The vector spaces
$\hat{V}_{f_i}$ consist of (0,1) forms  and together define a single
finite-dimensional vector space $\mbox{span}\{\hat{V}_{f_i}\}$. Let
${\mathbf V}$ denote the trivial vector bundle ${\mathbf
V}=\mbox{span}\{\hat{V}_{f_i}\}$ over $\N=\cup \N_{f_i}$. We then
have an exact sequence of locally trivial vector bundles over
$\N$\,:
$$
0 \ \longrightarrow\  \ker \sigma\   \longrightarrow\  {\mathbf V}\
\overset{\si}{\longrightarrow}\  \cok L \ \longrightarrow\  0
$$
where $\sigma$ is the composition of  the restriction map
(\ref{restriction}) with the projection to $\cok L$. Now, using a
metric on the vector bundle ${\mathbf V}$ induced from metrics on
$\CU$ and $N_D$, we define the obstruction bundle $\Ob$ over $\N$ to
be the orthogonal complement of $\ker\si$ in ${\mathbf V}$. Then
$\sigma$ restricts to a vector bundle isomorphism
$\Ob\overset{\cong}\to \cok L$ over $\N$. Lastly, the fact $\Ob$ has a canonical orientation is
proved in Lemma~\ref{orientable} below.
This completes the proof of Theorem~\ref{constuctionofOb}.

\vspace{1cm}

\setcounter{equation}{0}
\section{A generalized Image Localization lemma}
\label{section7}
\medskip

In this section we generalize the Image Localization Lemma of
Section~1, showing that it applies not just to solutions of
$\del_{J_\alpha}f=0$ but, more generally, to  maps $f:C\to N_D$ that
satisfy  of the perturbed $J_\alpha$-holomorphic map equation
$\del_{J_\alpha}f=\nu$ provided that $\nu$ is  small and its normal
component  lies in the fiber of the obstruction bundle.    The proof
uses a renormalization argument similar to those in Sections 6 and 7
of \cite{IP1}.  In the statement of the theorem, and throughout this
section, we will use the decomposition (\ref{split}) to write the
perturbation $\nu$  as the sum $\nu^{\ss T}+\nu^{\ss N}$ of tangent
and normal components.   The conclusion is the same as in the
original  localization lemma:  the images of the maps lie in the
divisor of $\alpha$, which is the zero section of the bundle $N \to
D$.

\medskip

\begin{theorem}[Image Localization with Perturbations]
\label{ImageinD} There is  a $\delta_0>0$ such that if
$\nu=\nu^T+\nu^N$ with   $|\nu| \leq \delta_0$ pointwise and
$\nu^N\in\Gamma(\Ob)$,  then  the image of every map $f:C\to N_D$
satisfying $\del_{J_\alpha}f=\nu$ lies in $D$.
\end{theorem}

\pf  If this statement is false, there is a sequence  $\{\nu_n\}$ of
the stated form with  $|\nu _n| \to 0$ pointwise,  and a sequence of
maps $f_n:C_n\to N_D$, each with at least one point not mapped into
$D$. By the compactness of the space of stable maps we may assume,
after passing to a subsequence, that the $f_n$ converge,  in
$C^0\cap W^{1,2}$ and in $C^\infty$ away from the nodes,  to a
$J_\alpha$-holomorphic map $f_0:C_0\to N_D$.  By Lemma~\ref{ILLemma}
the image of $f_0$  lies in  $D\subset N_D$.

We can now renormalize as in Lemma~6.3 of \cite{IP1}.      Let
$\phi_n$ be the projection of $f_n$ into $D$ and let $R_n:N_D\to
N_D$ be dilation in the fibers by a factor of $1/t_n$, where $t_n$
is chosen so that the normal component of the renormalized maps
$F_n=R_n\circ f_n$ has $C^1$ norm equal to 1.  Because the normal
component of $f_n$ converges to zero pointwise we have $t_n\to 0$,
and  since $\alpha$ vanishes along $D$, the pullbacks
$R_n^*K_\alpha$ converge to $0$ on compact sets in $N_D$.

Next write $F_n$  as the graph $(\phi_n, \xi_n)$ in $N_D$.  Then the
sections $\xi_n$ of $\phi^*_nN$ have $C^1$ norm 1 and  the original
maps are given by $f_n=(\phi_n, t_n\xi_n)$.  After trivializing the
pullback $\pi^*N$ by parallel transport along the fibers of
$\pi:N_D\to D$, we can compare the normal components of the (0,1)
forms $\Phi(f)=\del_{J}f-K_\alpha\bd f j$ for $f=f_n$ and $f=\phi_n$
on $C_n$  by
\begin{equation}
\label{lineardifference} \Phi^N(\phi_n)-\Phi^N(f_n) =
L_{f_n}(t_n\xi_n) + O(|t_n\xi_n|^2).
\end{equation}
But  $\Phi^N(\phi_n)=0$ since the image of $\phi$ lies in $D$ and by
assumption $\Phi^N(f_n)=\nu^N_{f_n}$ lies in the fiber of $\Ob_N$ at
$f_n$.   Dividing (\ref{lineardifference}) by $t_n$ and noting that
$\Ob_{f_n}$ is transverse to the image of $L_f$, we conclude that
\begin{equation}
\label{lineardifference2} \left|t_n^{-1}\nu^N_{f_n}\right| \to  0
\qquad \mbox{and}\qquad L_{f_n}(\xi_n) \to 0
\end{equation}
as $n\to\infty$.  But then the renormalized maps $F_n$ satisfy
$\Phi^N(F_n)=t_n^{-1}\nu^N_{f_n}\to 0$  and hence, as in the proof
of Lemma~6.3 of \cite{IP1}, converge to a limit $F_0:C_0\to N_D$ in
$C^0\cap W^{1,2}$ and in $C^\infty$ away from the nodes of $C_0$. If
$C_0$ is  smooth, we then have $f_n \to f_0$ and  $\xi_n\to \xi_0$
in $C^\infty$, and (\ref{lineardifference2}) implies  that $\xi_0$
is a non-trivial solution of
\begin{equation}
\label{L=0} L_{f_0}(\xi_0)  = 0.
\end{equation}
This contradicts the Vanishing Theorem~\ref{Vtheorem}.

When $C_0$ is nodal, we can restrict attention to one irreducible
component $C\subset C_0$ in which $\xi_0$ is non-trivial.  On each
compact subset of $C\ssetminus\{\mbox{nodes}\}$ we again have
$\xi_n\to \xi_0$ in $C^\infty$ and hence  $L\xi_0=0$ for
$L=L_{f_0}$.  After taking a sequence of such compact sets and
passing to a diagonal subsequence we obtain a non-trivial continuous
solution of $L\xi_0=0$ on $C\ssetminus\{\mbox{nodes}\}$.  We can then
see that $\xi_0$ is a weak solution of $L^*L\xi_0=0$ on $C$  as
follows.

Given  a smooth $\eta\in \Omega^{0,1}(C, \phi_0^*N)$, we must show
that the $L^2$ inner product $\langle L^*L\eta, \xi_0\rangle$
vanishes.  For each $\delta>0$, let $B(\delta)$ be the union of the
$\delta$-disks around the nodes of $C$ and let $\beta=\beta_\de$ be a cutoff function as defined before (\ref{defEbar}).
Writing $\xi_0=\beta\xi_0+(1-\beta)\xi_0$ and noting that
$|L((1-\beta)\xi_0)|\leq |d\beta||\xi_0|\lq c\de^{-1}|\xi|$,  we have the following bounds on the $L^2$ inner product:
\begin{eqnarray*}
|\langle L^*L\eta, \xi_0\rangle| &  \leq & |\langle L^*L\eta, \beta
\xi_0\rangle| +
|\langle L\eta, L(1-\beta)\xi_0\rangle| \\[.2cm]
& \leq& \big( \|L^*L\eta\|_{2; C} +c\delta^{-1}
\|L\eta\|_{2; B(\delta)} \big)\,  \|\xi_0\|_{\infty}\
\sqrt{\mbox{vol}(B(\delta))}.
\end{eqnarray*}
The right-hand side vanishes as $\delta\to 0$.  Thus $\xi_0$ is a
bounded weak solution of $L^*L\xi_0=0$ on $C$.

Finally, after a conformal change of metric, we can choose
coordinates around each node in which $L$ has the form $\del +A$
where $A$ is a zeroth-order operator and $L^*L=\Delta +B$  where $B$
is a first order operator.  Standard elliptic  theory then implies
that $\xi_0$ extends across the nodes to a smooth solution of
$L^*L\xi_0=0$ on $C$. Taking the inner product with $\xi_0$ and
integrating by parts then shows  that $\xi_0$ satisfies (\ref{L=0})
on $C$, which again contradicts   Theorem~\ref{Vtheorem}. \qed

\vspace{1cm}

\setcounter{equation}{0}
\section{The proof of the Main Theorem}
\label{section8}
\medskip

This section presents the proof of  the Main Theorem  stated in the
introduction: the local GW invariants of a spin curve $(D,N)$ are
given by the cap product
\begin{equation}\label{MainThm}
GW_{g,n}^{loc}(N_D,d)\,=\,\hev_*\big(\,[\bM_{g,n}(D,d)]^{\vir} \cap e(\Ob)\,\big)
\end{equation}
where $\Ob$ is the obstruction bundle  defined in Section~6.
The basic idea is to turn on a generic Ruan-Tian perturbation of the type described in
 Theorem~\ref{ImageinD}.  The perturbation defines a section  of $\Ob$ whose zero set represents
 the Euler class $e(\Ob)$, and on the other hand is cobordant to a cycle representing  the  local GW invariant.
The proof consists of three steps.

\bigskip

To simplify notation, we set
\begin{equation*}
\begin{array}{ll}
\Map_D \,=\, \Map_{g,n}(D,d),\ \ \ \ \ \ \  & Y_D\,=\,\CM_{g,n}\times D^n_{\textstyle\phantom{\sum}} \\
\ \ \Map \,=\, \Map_{g,n}(N_D,d),    & \ \, Y\,=\, \CM_{g,n}\times {N_D^n}
\end{array}
\end{equation*}
Let $\hat{\F}\to \Map$ be the topological vector bundle
whose fiber at $f$ is $\Omega^{0,1}(f^*TN_D)$. This bundle has a
section $\Phi$ given by $J_\a$-holomorphic map equation
$\Phi(f)=\del_J f - K_\a( \bd_J f ) j.$ For each perturbation $\nu$
in ${\cal P}$ defined in (\ref{defP}), one can perturb  $\Phi$ to
obtain a section
\begin{equation}
\label{8.defPhi}
 \Phi_\nu(f)\,=\,\Phi(f)-\nu_f
\end{equation}
where $\nu_f$ is the restriction of $\nu$ to the graph of $f$ as in
(\ref{restriction}). The zero set $\CM^\nu=\Phi_\nu^{-1}(0)$ is the
moduli space of $(J_\a,\nu)$-holomorphic maps into $N_D$.

By Theorem~\ref{constuctionofOb}, the bundle $\Ob$ is defined on a  neighborhood $\N$ of the moduli space of $J_\a$-holomorphic maps inside  $\Map$.
Following \cite{RT2}, one can decompose $\N$ into strata $\N_B$
indexed by a finite collection of sets $B$ that specify: (i) a homeomorphism type of the
domain as a curve with marked points and (ii) a degree of map
associated to each component of the domain.

\bigskip\medskip

\noindent{\bf Step 1.\,}
As in  Section~7, we can decompose
$\nu$   according to the splitting of
$f^*TN_D$ into tangent and normal subspaces.
Fix a small perturbation $\nu_{\ss T}$ that is generic in the subspace $\CP_T\subset\CP$ of perturbations that lie in the tangential subspace and consider the section $\Phi_{\nu_{\ss T}}$ defined by (\ref{8.defPhi}).

By the Image Localization
Theorem~\ref{ImageinD} the images of all $(J_\a,\nu_{\ss T})$-holomorphic maps lie in $D$.
The transversality arguments of
Section 3 in \cite{RT2}, applied to maps into  $D$,
imply that the moduli space
$$\CM\,=\,\Phi_{\nu_{\ss T}}^{-1}(0)\,\subset\,\N\cap\Map_D
$$
has a stratification such that each stratum
$\M_B=\CM\,\cap\, \N_B$
is a smooth oriented manifold of dimension $4\b + 2n -2k_B$ where
$\b=d(1-h)+g-1$ as in (\ref{defbeta}) and $k_B$ is the number of
nodes of the domains of maps in $\N_B$.

Let $\M$ be the top stratum consisting of maps with smooth domains.  Now apply the construction described at the end of Section~3:
fix a neighborhood $U$ of the image $\hev(\bd\CM)$ of the boundary $\bd\CM=\CM\ssetminus \M$ and consider the manifold with boundary
$\bM_U\subset \M$ as in Definition~\ref{defn3.1}.   The following lemma shows that space $\bM_U$ represents  the Li-Tian  virtual fundamental class (\ref{LTVFC})  modulo  small neighborhoods of its boundary.

\begin{lemma}\label{rel-VFC}
Let $\hat{U}=\hat{ev}^{-1}(\overline{U})$. Then
$$
j_*[\CM_{g,n}(D,d)]^{\vir}\,=\,i_*[\CM_U]\,\in\, H_{4\b+2n}(\Map_D,\hat{U})
$$
where $j: \Map_D\to (\Map_D,\hat{U})$ and $i:(\CM_U,\bd\CM_U)\to (\Map_D,\hat{U})$
are inclusion maps.
\end{lemma}

\pf Let $\hat{\F}_D\to \Map_D$ be the bundle whose fiber at $f$ is $\Omega^{0,1}(f^*TD)$.
Since $K_\a\equiv 0$ on $D$, we can regard $\Phi_{\nu_{\ss T}}$ as a section of $\hat{\F}_D$.
Proposition 3.4  of \cite{LT1}  then shows how
the section $\Phi_{\nu_{\ss T}}$ gives rise to
a cover of  the moduli space $\CM=\Phi_{\nu_{\ss T}}^{-1}(0)$ by  finitely many {\em smooth approximations}
$\{(W_k, F_k)\}$.  This means that
\begin{itemize}
\item
each $W_k$ is open in $\Map_D$ and $\CM\subset
\bigcup W_k$, and
\item
each $F_k$ is a subbundle of $\hat{\F}_D$ over $W_k$ with finite rank such
that $\Phi_{\nu_{\ss T}}^{-1}(F_k)\subset W_k$ is smooth and $F_k$ restricts to a
smooth bundle over $\Phi_{\nu_{\ss T}}^{-1}(F_k)$ with smooth section ${\Phi}_{\nu_{\ss T}}$.
\end{itemize}
Observe that, because the top stratum $\M$ is already smooth, we can assume
\begin{itemize}
\item $W_1$ consists of all maps with smooth domains (so $W_1\cap\bM=\M$) and $F_1$ has rank zero.
\end{itemize}

The proofs of  Proposition 2.2 and Theorem 1.2 of \cite{LT1}
describe how one can perturb the moduli space $\CM$
inside the union of the smooth manifolds $\Phi^{-1}(F_k)$ to obtain
a cycle that represents the virtual fundamental class
$[\CM_{g,n}(D,d)]^{\vir}$.

Now choose  a neighborhood $V$ of  the image $\hat{ev}(\bd\CM))$ in $Y_D$ with $\overline{V}\subset U$.  Set
$\hat{V}=\hev^{-1}(V)$ and  replace each $W_i$ with $i\geq 2$ with $W_i\cap\hat{V}$.  Then
the open set $W= \bigcup_{k>1} W_k$ lies in $\hat{V}$.
Following \cite{LT1} we can then perturb  $\CM \cap W$ while
keeping $\bM$ fixed outside of $\hat{V}$ -- hence on $\bM_U$ -- to produce a cycle
representing $[\CM_{g,n}(D,A)]^{\vir}$.  Passing to homology proves the statement of the lemma.
\qed

\begin{rem} \label{RT=LT}
The virtual fundamental class that appears in Lemma~\ref{rel-VFC} is the Li-Tian class (\ref{LTVFC}).  We can verify that it is also the Ruan-Tian class:
Let $k:Y_D\to (Y_D,\overline{U})$ be an inclusion map.
Since $\hev\circ j=k\circ \hev$, Lemma~\ref{rel-VFC} gives $k_*\circ \hev_*[\CM_{g,n}(D,d)]^{\vir}=\hev_*[\CM_U]$.
The uniqueness statement in Lemma~\ref{KMGW} then shows that
$$
\hev_*[\CM_{g,n}(D,d)]^{\vir}\,=\,GW_{g,n}(D,d)\,\in\,H_{4\b+2n}(\CM_{g,n}\times D^n;\Q).
$$
\end{rem}

\bigskip\bigskip

\noindent{\bf Step 2.\,} We now further perturb the section
$\Phi_{\nu_{\ss T}}$ by adding a section of the obstruction bundle
$\Ob$ induced from a Ruan-Tian perturbation $\mu$.   Each $\mu\in{\cal P}$ defines a section $\hat{s}_\mu$ of the
obstruction bundle $\Ob$ over $\N$ by
$$
\hat{s}_\mu(f)\,=\,P_f(\mu_f)
$$
where $P_f$ is the $L^2$-orthogonal projection onto the fiber
$\Ob_f$.   Let $s_\mu$ denote the restriction of $\hat{s}_\mu$
to the moduli space $\CM$. The Image Localization
Theorem~\ref{ImageinD} implies that
\begin{equation}
\label{zero-set} (\Phi_{\nu_{\ss T}} - \hat{s}_\mu)^{-1}(0)\,=\,
\CM\cap \hat{s}_\mu^{-1}(0)\,=\, s_\mu^{-1}(0).
\end{equation}

\begin{lemma}\label{Trans-1}
For generic $\mu$ in ${\cal P}$, the space (\ref{zero-set})  has a
stratification indexed by  $B$ such that each stratum
$\M_B\cap s_\mu^{-1}(0)$ is a smooth manifold of
dimension $2\b+2n-2k_B$.
\end{lemma}

\pf The proof is a transversality argument  using the universal
moduli space over the space of perturbations ${\cal P}$ (cf.
Theorem~3.1 of \cite{RT2}). The smooth bundle
$$
\Ob\,\to\,  \M_B\times {\cal P}
$$
has rank $2\beta$ and has a section $\sigma_B$ defined by
$\sigma_B(f,\mu)=s_\mu(f)$  whose zero set is the universal moduli
space associated with $B$.  The differential of $\sigma_B$  at $(f,\mu)$ is
$$
(D\sigma_B)_{(f,\mu)}(\xi,\chi)\,=\,(Ds_\mu)_f(\xi)- P_f(\chi_f)\ \
\ \ \ \ \mbox{for }\ \  \xi\in T_f\M_B,\
 \chi\in T_\mu{\cal P}
$$
where $Ds_\mu$ is the differential of the restriction of $s_\mu$ to
$\M_B$ and $\chi_f$ is the restriction of $\chi$ to
the graph of $f$ as in (\ref{restriction}). Since the map $\chi\to
P_f(\chi_f)$ is onto, so is the differential
$(D\sigma_B)_{(f,\mu)}$. The universal moduli space ${\cal
U}=\sigma_B^{-1}(0)$ is therefore smooth. Now, consider the projection
$$
\pi_B: {\cal U} \to {\cal P}
$$
given by $\pi_B(f,\mu)=\mu$. Since $\pi_B$ is Fredholm, the
Sard-Smale Theorem  implies that for generic $\mu$ the fiber
$\pi_B^{-1}(\mu)=\M_B\cap s_\mu^{-1}(0)$ is a smooth manifold of dimension
$\dim\!\M_B - \mbox{rank }\Ob $, which is $2\b + 2n -2k_B$. \qed

\medskip

Fix a generic $\mu\in \CP$ as in Lemma~\ref{Trans-1} and set
\begin{equation}\label{fixed-section}
\bar{Z}\,=\,s_\mu^{-1}(0)\,=\,(\Phi_{\nu_{\ss T}} - \hat{s}_\mu)^{-1}(0).
\end{equation}
By Lemma~\ref{Trans-1} and the discussion  above Lemma~\ref{KMGW}, applied to some neighborhood of
$\hev(\bd \bar{Z})$ where $\bd \bar{Z}=\bar{Z}\cap \bd \CM$,
the image $\hat{ev}(\overline{Z})$ defines a homology class
\begin{equation}\label{zero-class}
[\,\hev(\bar{Z})\,]\,\in\,H_{2\b+2n}(Y_D;\Q)\,\cong\,H_{2\b+2n}(Y;\Q).
\end{equation}
On the other hand, the local invariant GW class is determined by the
image of the  moduli space $\CM^\nu=\Phi_{\nu}^{-1}(0)$
for a generic $\nu\in {\cal P}$ as in (\ref{LGWclass})\,:
\begin{equation}\label{localGWclass}
GW_{g,n}^{loc}(N_D,d) \,=\,    \big[\hev(\CM^\nu)\big]
\in\,H_{2\b+2n}(Y;\Q).
\end{equation}
Thus the local invariant is defined in terms of the zero set of
$\Phi_\nu$ for a fixed generic  $\nu\in {\cal P}$, while the
class (\ref{zero-class}) is defined in terms of the zero set of
$\Phi_{\nu_{\ss T}}-\hat{s}_\mu$.   To show these are equal we introduce, for
each pair $(\nu,\mu)$ in ${\cal P}\times{\cal P}$, a section of the
bundle $\hat{\E}\to\Map$ defined by
\begin{equation*}
\Phi_{\nu,\upsilon}(f)\, =\, \Phi(f) -\nu_f - \hat{s}_\upsilon(f).
\end{equation*}
The transversality argument of Lemma~\ref{Trans-1},   now applied to the
universal moduli space over the parameter space ${\cal P}\times{\cal
P}$, shows that,  for generic small $(\nu,\upsilon)$ in ${\cal P}\times{\cal P}$,
$\Phi_{\nu,\upsilon}^{-1}(0)$ has a stratification indexed by $B$ such
that each stratum $\N_B\cap \Phi_{\nu,\upsilon}^{-1}(0)$ is a smooth
manifold of dimension $2\b+2n-2k_B$. Standard cobordism argument (cf. Theorem 3.3 of \cite{RT2}) then give
\begin{equation}\label{ThmA}
GW_{g,n}^{loc}(N_D,d)\,=\,[\,\hev (\Phi_{\nu,0}^{-1}(0))\,]\,=\,
[\,\hev(\Phi_{\nu_{\ss T},\mu}^{-1}(0))\,]\,=\,[\hev(\bar{Z})].
\end{equation}

\bigskip\bigskip

\noindent{\bf Step 3.\,}
It remains to show that $[\hev(\bar{Z})]$ equals to the right-hand side of (\ref{MainThm}).
Again, to avoid issues of  smoothness near $\partial\bM$, we work in relative homology.
The following lemma is, in essence,  the statement  that the zero set of a generic section of a vector bundle
is Poincar\'{e} dual (in relative homology) to the Euler class.

\begin{lemma}
Let $\bar{Z}_U=\CM_U\cap \bar{Z}$. Then
\begin{equation}\label{PD-key}
j_*(\,[\CM_{g,n}(D,d)]^{\vir}\cap e(\Ob)\,)\,=\,i_*[\bar{Z}_U]\,\in\,H_{2\b+2n}(\Map_D,\hat{U}).
\end{equation}
where $j: \Map_D\to (\Map_D,\hat{U})$ and $i:(\M_U,\bd\M_U)\to (\Map_D,\hat{U})$
are inclusion maps.
\end{lemma}

\pf The obstruction bundle $\Ob$ restricts to a smooth bundle $\Ob_U$
over $\CM_U$.
Let $s_{\ss U}$ be the restriction of the section $\hat{s}_\mu$ in
(\ref{fixed-section}) to $\CM_U$.
Without loss of generality, we can assume
$s_{\ss U}$ is transverse to the zero section, so $\bar{Z}_U=s^{-1}_{\ss U}
(0)$ is a compact smooth manifold with boundary.

Define a double $\tilde{\M}_U$ of $\CM_U$ by identifying two copies of $
\CM_U$ along $\bd \CM_U$.
Similarly, define
a vector bundle $\tilde{\Ob}_U\to \tilde{\M}_U$ and its section $
\tilde{s}_{\ss U}$ that are doubles of
$\Ob_U$  and $s_{\ss U}$ respectively. Since $\tilde{\M}_U$ is a closed
manifold, we have
\begin{equation}\label{PD-2}
[\tilde{\M}_U]\,\cap\,e(\tilde{\O}b_U)\,=\,
[\tilde{Z}_U]
\end{equation}
where $\tilde{Z}_U=\tilde{s}_{\ss U}^{-1}(0)$. Consider the following
commutative diagram\,:
\begin{equation}\label{PD-3}
\xymatrix{
H_{4\b}(\tilde{\M}_U)\ \     \ar[r]^{\kappa_*\ \ \ \ \ \ \ \ }  \ar[d]^{\,
\cap\, e(\tilde{\Ob}_U)} &
\ \ H_{4\b}(\tilde{\M}_U,\tilde{\M}_U\ssetminus \M_U^\circ)\ \
\ar[d]^{\, \cap\, e(\tilde{\Ob}_U)} &
\ \ H_{4\b}(\M_U,\bd \M_U)  \ar[l]_{\ \ \ \ \   \iota_*}^{\ \ \ \ \
\simeq} \ar[d]^{\,\cap\,e(\Ob_U)} \\
H_{2\b}(\tilde{\M}_U)\ \  \ar[r]^{\kappa_*\ \ \ \ \ \ \ \ } &
\ \ H_{2\b}(\tilde{\M}_U,\tilde{\M}_U\ssetminus \M_U^\circ)\ \  &
\ \ H_{2\b}(\M_U,\bd \M_U)  \ar[l]_{\ \ \ \ \   \iota_*}^{\ \ \ \ \
\simeq} }
\end{equation}
where $\M_U^\circ=\M_U\ssetminus \bd\M_U$, $\kappa$ and $\iota$ are inclusion maps,
each rectangle commutes by the naturality of
the cap product, and the horizontal isomorphisms follow by excision.
Consequently, (\ref{PD-2}), (\ref{PD-3}) and the
facts
$\kappa_*[\tilde{\M}_U]\,=\,\iota_*[\CM_U]$ and $\kappa_*[\tilde{Z}_U]\,=\,\iota_*[\bar{Z}_U]$ show
\begin{equation*}
[\CM_U]\,\cap\,e(\Ob_U)\,=\,[\bar{Z}_U]\,\in\,H_{2\b+2n}(\M_U,\bd\M_U).
\end{equation*}
Now, (\ref{PD-key}) follows from this equality, Lemma~\ref{rel-VFC} and the naturality of cap product.
\qed

\medskip

\bigskip
\non
{\bf Proof of Main Theorem.\,}
Recall that $\CM$ and $\bar{Z}$ are decomposed into smooth strata indexed by sets $B$ such that
$\dim\!\M_B=4\b+2n-2k_B$ and $\dim\!Z_B=2\b+2n-2k_B$ where $Z_B=\bar{Z}\cap \M_B$.
Because $Y_D$ is a compact manifold, any class in $H_*(Y_D;\Q)$ can be represented by an embedded submanifold by Thom's Theorem \cite{Th}.
We can thus choose a basis for $H_{m-2\b}(Y_D;\Q)$, where $m+2n=\dim\!Y_D$,  represented by submanifolds $D_i$  in general position with respect to all the restriction maps $\hev_{|\M_B}$ and $\hev_{|Z_B}$.   Counting dimensions, one sees that
\begin{itemize}
\item each $D_i$ is disjoint from the image $\hev(\bd \bar{Z})\subset \hev(\bd \CM)$, and
\item we can choose   a submanifold $D$ representing the class $\hev_*([\CM_{g,n}(D,d)]^{\vir}\cap e(\Ob))$ that   is
disjoint from $D_i\cap \hat{ev}(\bd\CM)$ for all $i$.
\end{itemize}
Hence, by shrinking the neighborhood $U$ of $\hev(\bd\CM)$  if necessary, we can assume that

\begin{equation}\label{key-int}\tag{8.11ab}
\addtocounter{equation}{1}
\overline{U}\cap \hat{ev}(\bar{Z})\cap D_i\,=\,\emptyset\ \ \ \ \ \ \
\mbox{and}\ \ \ \ \ \ \
\overline{U}\cap D\cap D_i=\emptyset
\end{equation}
for all $i$.  Now, consider the commutative diagram\,:
\begin{equation}\label{key-diag}
\xymatrix{
 H_{2\b}(\Map_D)  \ar[r]^{j_*\ \ } \ar[d]^{\hev_*} &
 H_{2\b}(\Map_D,\hat{U})  \ar[d]^{\hev_*} &
 H_{2\b}(\CM_U,\bd\CM_U)   \ar[l]_{\ \ \ i_*}  \\
 H_{2\b}(Y_D)  \ar[r]^{k_*\ \ }  &
 H_{2\b}(Y_D,\bar{U}) & }
\end{equation}
Observe that (\ref{key-int}) implies that for all $i$
\begin{equation}\label{8.11}
(a)\ [D]\cdot [D_i]\,=\,k_*[D]\cdot k_*[D_i]
\ \ \ \ \ \mbox{and}\ \ \ \ \
(b)\ [\hev(\bar{Z})]\cdot [D_i]\,=\,\hat{ev}_*\circ i_*([\bar{Z}_U]) \cdot k_*[D_i]
\end{equation}
where, in both of these equations,  the dot on the right-hand side  of is
the intersection pairing in $H_{2\b}(Y_D,\bar{U})\cong H_{2\b}(Y_D/U, \bd U)$.   By the definition of $D$, (\ref{8.11}a) states that
$$
\hat{ev}_*(\,[\CM_{g,n}(D,d)]^{\vir}\cap e(\Ob)\,)\cdot [D_i]
\ =\ k_*\circ \hat{ev}_*(\,[\CM_{g,n}(D,d)]^{\vir}\cap e(\Ob)\,)\cdot k_*[D_i]
$$
But the intersection on the right is, by the commutativity of the diagram,  (\ref{PD-key}) and (\ref{8.11}b),
$$
\hat{ev}_*\circ j_*(\,[\CM_{g,n}(D,d)]^{\vir}\cap e(\Ob)\,)\cdot k_*[D_i]
\ =\
\hat{ev}_*\circ i_*([\bar{Z}_U]) \cdot k_*[D_i]\,=\,
[\hev(\bar{Z})]\cdot [D_i].
$$
This shows $\hat{ev}_*(\,[\CM_{g,n}(D,d)]^{\vir}\cap e(\Ob)\,)=[\hev(\bar{Z})]$ and hence,
together with  (\ref{ThmA}), completes the proof of the Main Theorem (\ref{MainThm}). \qed

\vspace{1cm}

\setcounter{equation}{0}
\section{Secondary index invariants}
\label{section9}
\medskip

This section puts the obstruction bundle in a general context and
explain why Euler class $e(\Ob)$  cannot in general be computed
using Grothendieck-Riemann-Roch or the Families Index theorems.  It
also includes the calculation of the Euler class in some examples.

Let $E$ and $F$ be Banach vector bundles over a compact parameter
space $X$.  We can consider the vector bundle
$$
\begin{CD}
\mbox{Fred}_\ell(E,F)\\
@VVV\\
X
\end{CD}
$$
whose fiber over $x\in X$ is the space of  real linear Fredholm maps
from $E_x$ to $F_x$ with index $-\ell$.  A  section $L$ of
$\mbox{Fred}_\ell(E,F)$ then defines a family-index class
$$
\mbox{ind}\, L \in KR(X).
$$
obtained by pulling back the class of the virtual bundle $[\ker
L]-[\cok L]$ on $\mbox{Fred}_\ell(E,F)$.  The index theorem for
families \cite{AS} gives formulas for the Pontryagin classes
$p_i(\mbox{ind}\, L)$.  But the Euler class does not factor through
K-theory, and hence $e(\mbox{ind}\, L)$  is not computable in the
same way and in fact is not even defined  in general.

On the other hand,  for each $k>0$, the set $A_{k, \ell}$ of all
$L\in \mbox{Fred}_\ell(E,F)$ with $\dim \ker L= k$ and $\dim \cok L
= k+\ell$ is a submanifold of  real codimension $k(k+\ell)$.  As
shown in \cite{K},  the closures $\overline{A}_{k, \ell}$ of these
submanifolds define ``Koschorke classes''
\begin{equation}
\label{KoschorkeClasses} \kappa_{k,\ell} \in
H^{k(k+\ell)}(\mbox{Fred}_\ell(E,F)).
\end{equation}
A  section $L$ of $\mbox{Fred}_\ell(E,F)$ then defines classes
$L^*\kappa_{k,\ell} \in H^{k(k+\ell)}(X)$.  The classes
$\{L^*\kappa_{k+i,\ell}\, |\, i>0\}$  are the obstructions to
deforming $L$ within its homotopy class to a family $\{L_x\}$ of
operators with $\dim \ker L_x \leq k$ for all $x\in X$.   In
particular, when
all  Koschorke classes  $\{L^*\kappa_{k,\ell}\}$ vanish,  $L$ can be deformed  to  a section $\tilde{L}$ with $ \ker \tilde{L}_x =0$ for all $x\in X$; $\cok \tilde{L}$ is then a rank $\ell$ vector bundle over $X$ that represents $\mbox{ind}\, L \in KR(X)$. In this sense, the  Koschorke classes are the obstruction to realizing the family index --- which is defined as a formal difference of bundles ---  as an {\em actual} vector bundle.   Furthermore, when the Koschorke classes   vanish this bundle is well-defined up to homotopy, so the Euler class
$$
e(\mbox{ind}\, L) \in H^\ell(X)
$$
is defined. This is a ``secondary class'' in the sense that it
exists only for families with vanishing Koschorke classes.

Now consider the situation at hand, where
$X=\bM_{g,n}(D,d)\subset \Map_{g,n}(N_D,d)$, the fibers of $E$ and $F$ are the Sobolev completions of $\Omega^0(f^*N)$ and $\Omega^{0,1}(f^*N)$ and $L$ is the linearization map $f\mapsto L_f$ with index $\ell=2\beta$.   Then

\begin{itemize}
\item  $L^*\kappa_{k, 2\beta}=0$ for all $k$  by the Vanishing Theorem~\ref{Vtheorem}.
\end{itemize}
This gives a global perspective on the role of the  Vanishing
Theorem~\ref{Vtheorem}:  it ensures that all  Koschorke classes
vanish.  In fact, it shows that after perturbing the K\"{a}hler
structure $J$ to $J_\alpha$,  the space of $J_\alpha$-holomorphic
maps is mapped by $\Psi(f)=L_f$ into a region in
$\mbox{Fred}_\ell(E,F)$ where the index bundle is an actual bundle.
The pullback $\Psi^*(\mbox{ind}\, L)$ is the obstruction bundle
$\Ob$ of the Main Theorem.

\bigskip

Our proof of the Main Theorem in Section~8 used the fact that the
obstruction bundle carries a canonical  orientation.  The proof,
which we give now, fits into the above discussion of the space of
Fredholm operators.
Below, we write  $\CM$ for $\CM_{g,n}(D,d)$.

 \begin{lemma}
\label{orientable} The bundle $\Ob$ is  orientable and has a
canonical  orientation.
\end{lemma}

 \pf    The linearization $f\mapsto L_f$ defines a map $L:\bM\to \mbox{Fred}$ whose image lies in the set of Fredholm operators with trivial kernel.  Over this set,  $\Ob$ -- whose fiber is the cokernel ---   is isomorphic to the dual of the index bundle.   The obstruction to orientability is therefore the pullback   $w=L^*w_1(\mbox{det}_\R\  L)\in H^2( \bM,\Z_2)$ of the  first Stiefel-Whitney class of the real determinant of the  index bundle.  But the map $L$ extends canonically to a homotopy $\bM\times [0, 1] \to \mbox{Fred}$ by writing $L_f= \del+R_\alpha$ and setting $L_t(f)=  \del+t R_\alpha$.  Because $\mbox{det}_\R\  L$ is defined over the entire space of Fredholm operators, we have $w=L_0^*w_1(\mbox{det}_\R\  L)$.  This is zero because the image of $L_0$ consists of complex Fredholm operators, whose kernel and cokernels have canonical complex orientations.

This shows, in fact, that the orientation bundle
$\Lambda^{top}\Ob$ is trivial over $\bM\times [0, 1]$.  There are
therefore two orientation classes (nowhere-vanishing  sections
modulo multiplication by positive functions).  The one that agrees
with the complex orientation along $\bM\times \{0\}$ will be called
the ``canonical''  orientation on $\Ob$. \qed

\bigskip

Finally, we must specify the orientation on $\Ob$ for which
the Main Theorem holds.  Let $\bar{Z}$ be the zero set of a
transverse section $\Ob \to \bM$.  At each $(C,f)\in \bar{Z}$, the standard orientations of $\CM$ and $\bar{Z}$ used in GW theory 
are given by the determinant bundles $\det (\del^T\oplus Jdf)$ and $\det (\del^T\oplus \del\oplus Jdf)$
respectively. Thus in the equality (\ref{lp3}) in the Main Theorem,  the cycles representing the two sides are
consistently oriented provided
\begin{equation*}
\det(\del^T\oplus Jdf)  \,=\,  \det(\del^T\oplus \del\oplus Jdf)  \otimes \Lambda^{top}\Ob
\end{equation*}
This equality holds since the canonical  orientation defined in the proof of
Lemma~\ref{orientable} is 
$$
\Lambda^{top}\Ob \,=\,\det(\del)^*.
$$

\noindent{\bf Remarks on calculations} \\

We conclude with some remarks on calculating the  Euler class of the
obstruction bundle.  Algebraic geometers have a  standard procedure
for calculating the Euler class of the index bundles of families of
$\del$ operators by using the Grothendieck-Riemann-Roch formula,
often in conjunction with localization by a group action. This
procedure is not, in general, applicable to finding the Euler class
of the {\em real} bundle $\Ob$.  But it is worth noting that the GRR
formula yields some information in the following two cases.

\medskip

\noindent (1)  Since the square of the Euler class is the top
Pontryagin class, we have $e^2(\Ob) = p_{\beta}(\Ob) = (-1)^\beta
c_{2\beta}(\Ob\otimes_\R \cx)$.  But $\Ob\otimes_R \cx$ is the
complex index bundle $\mbox{ind}_\cx L$ because $\ker L$ vanishes
for all operators in the family.  Thus the GRR formula is applicable
for finding  $e^2(\Ob)$, but not for finding $e(\Ob)$ itself.

\medskip

\noindent (2)  Recall Givental's notion of {\em twisted GW invariants}:
each bundle $E$ over $X$, determines a bundle ${\cal E}\to
\bM_{g,n}(X,A)$ over the space  of stable maps whose fiber over a
map $f:C\to X$ is $H^0(C,f^*E)\ominus H^{1}(C,f^*E)$.  The
corresponding twisted invariants are obtained by evaluating Chern
classes of ${\cal E}$, together with $\tau$ classes, on  the virtual
fundamental class. In some especially simple situations,  the
Euler class of the obstruction bundle can be expressed in terms of
twisted GW invariants, and these can be
calculated.

\begin{lemma}
\label{exactsequenceLemma} Over the space of  $J$-holomorphic maps
there is  a locally trivial complex vector bundle $Q$ that is
isomorphic to the obstruction bundle $\Ob$ as an {\em unoriented}
real vector bundle,  such that for each $J$-holomorphic map $f$ the
fiber $Q_f$  of $Q$ fits into the split exact sequence of complex
vector spaces
\begin{equation}
\label{exactsequence} 0 \longrightarrow \overline{\mbox{\ker}
\del_f} \overset{R_\a}\longrightarrow\mbox{\cok} \del_f
\longrightarrow Q_f \longrightarrow 0.
\end{equation}
\end{lemma}

\pf   We will repeated use the facts that $R_\a$ satisfies  $JR_\a=-R_\a J$ and
\begin{equation}\label{crossterms}
\langle \del\xi_1,R_\a\xi_2\rangle\,=\,-\langle R_\a\xi_1,\del\xi_2\rangle
\end{equation}
where $\langle\cdot,\cdot\rangle$ denote the $L^2$ inner product (cf. section 8 of \cite{LP}).
Fix a $J$-holomorphic map
$f:C\to D$ and write $\del_f$ as $\del$. Regard $\cok \del$
as the $L^2$ orthogonal complement  to the image of $\del$.    First
note that if $\xi\in \ker\del$ then by (\ref{crossterms})
$$
\langle  \del\xi', R_a\xi\rangle \ =\ - \langle  R_\a\xi',\del
\xi\rangle \ =\ 0
$$
for all $\xi'\in \Omega^0(f^*TD)$.  Hence $R_\a\xi$ lies in $\cok
\del$.  Defining $Q$ to be the $L^2$ perpendicular
$$
Q\ =\ \left(R_a(\ker\del)\right)^\perp \, \subset\, \cok \del
$$
gives the exact sequence (\ref{exactsequence}) of complex  vector
spaces.  Note that $R_\a$, originally complex anti-linear, becomes
complex linear when we reverse the complex structure on $\ker \del$
and that this definition of $Q$ splits the sequence.

Now let $q:\cok (\del+R_a)\to \cok \del$ be the $L^2$ orthogonal
projection.   Observe:
\begin{itemize}
\item  $q$ is injective:  any $\eta\in\cok \del$ with $q(\eta)=0$ has the form $\eta=\del\xi$ for some $\xi$, and hence vanishes because,  by (\ref{crossterms}),
$$
0\, =\, \langle  \del\xi, (\del+R_a)\xi\rangle \,=\,
\|\del\xi\|^2\, =\, \|\eta\|^2.
$$
\item  The image of $q$ lies in $Q$:  for any $\eta\in\cok (\del+R_\a)$ and $\xi\in \ker \del$ we have
$$
0\, =\, \langle  \eta, (\del+R_a)\xi\rangle \, =\,  \langle  \eta,
R_a\xi\rangle \,=\,   \langle  q(\eta), R_a\xi\rangle
$$
(the last equality holds because $R_\a\xi\in\cok \del$ as above).
Thus $q(\eta)$ is $L^2$ perpendicular to the image of $R_\a$.
\end{itemize}
Finally, count dimensions.  From (\ref{exactsequence}) we have
$\dim Q_f= -\ind \del = -\ind (\del+R_\a)$.   But $\ker
(\del+R_\a)=0$ by Theorem~\ref{Vtheorem}, so $\dim Q_f= \dim
\cok (\del+R_\a)$.  Thus $q$ is an isomorphism onto its image $Q$. \qed

\bigskip

As $f$ varies across the space  $\CM=\CM_{g,n}(D,d)$ of
$J$-holomorphic maps,  one obtains families
$$\ker\del\,\to\,\CM\ \ \ \ \ \ \
\mbox{and}\ \ \ \ \ \ \ \cok\del\,\to\,\CM
$$
whose fibers are the complex vector spaces $\ker\del_f=H^0(f^*N)$
and $\cok\del_f=H^1(f^*N)$.   In general,  the dimensions of these
fibers is not  constant:  the dimension of the kernel and the
cokernel  jumps up (by equal amounts) along a ``jumping locus'' in
$\CM$.   But away from the jumping locus
Lemma~\ref{exactsequenceLemma} gives a formula, due to Kiem and Li
{\cite{KL},  for the Euler class of the obstruction bundle:

\begin{prop}\label{con-rank}
Suppose $\ker\del$ and $\cok\del$ are locally trivial vector
bundles over  a  set $Z\subset \CM$.  Then there is an isomorphism
of oriented real vector bundles
\begin{equation}
\label{orientedisom} \Ob \simeq (-1)^{h_i}\, Q
\end{equation}
over each component $Z_i$ of $Z$ where $h_i$ is $h^0(f^*N)$ on
$Z_i$. Consequently, $e(\Ob)\in H^*(Z)$ is
\begin{equation}
\label{constantranke(Ob)} e(\Ob)\,=\,\sum (-1)^{h_i}\,
c_{top}(\,[\cok\del_{\,\,|\,Z_i}]-[\overline{\ker\del}_{\,\,|\,Z_i}]\,).
\end{equation}
where the sum is over all connected components $Z_i$ of $Z$.
\end{prop}

\pf $Q$ has a complex orientation from (\ref{exactsequence}) and the
orientation of $\cok(\del+R_\a)=\Ob$ is given by $\det(\del)^*$ as
in Lemma~\ref {orientable}.  These are related by
\begin{align}\label{E-C-2}
\det(\del)^*\,
&=\,\wedge^{top}\cok\del\otimes \wedge^{top}\ker\del\notag \\
& \,=\, (-1)^{h^0} \wedge^{top}\cok\del\otimes \wedge^{top} \overline{\ker\del} \notag \\
&=\,(-1)^{h^0}\det Q,
\end{align}
which gives (\ref{orientedisom}).  Taking Euler classes and noting
that $Q$ is a complex bundle, we have   $e(\Ob)=(-1)^{h_i}
c_{top}(Q)$ and   (\ref{constantranke(Ob)}) follows from the exact
sequence of Lemma~\ref{exactsequenceLemma}. \qed

\medskip

\begin{ex}  Suppose that $D$ is an elliptic curve with odd theta characteristic.  Because the local GW invariants depend only on the parity of the theta characteristic, we can take $N$ to be a trivial bundle.  Then
$\ker\del$ is the trivial line bundle $\cx$ over $\CM$. In this
case, we have
\begin{equation}\label{h=1case}
e(\Ob)\,=\,-c_{top}(\,[\cok\del]-
[\overline{\cx}]\,)\,=\,-c_{top}(\cok\del)\,=\,c_{top}(\ind(\del)).
\end{equation}
\end{ex}
\medskip

The formula (\ref{h=1case}) shows that, in this case, the local GW invariants  are
special cases of Giventhal's twisted GW invariants of curves. They
can be explicitly computed using the result of Proposition~2  in the
paper of Faber and Pandharipande \cite{FP}.

\vspace{1cm}

{\small

\vspace{1.3cm}

\noindent {\em  Department of  Mathematics,  University of Central Florida, Orlando, FL 32816\\[.15cm]
Department of  Mathematics,  Michigan State University, East Lansing, MI 48824}

\medskip

\noindent {\em e-mail:}\ \ {\ttfamily junlee@mail.ucf.edu, \ \ \ {\ttfamily parker\@@math.msu.edu}

}

\end{document}